\documentclass{article}

\usepackage{hyperref}
\hypersetup{
    colorlinks=true,
    linkcolor=blue,
    filecolor=magenta,      
    urlcolor=cyan,
    pdftitle={Admissibility of mean bounds},
    pdfpagemode=FullScreen,
    }

\usepackage{soul}
\usepackage[utf8]{inputenc}
\usepackage{xcolor}
\usepackage{graphicx}
\usepackage{subcaption}
\usepackage{algorithm}
\usepackage{algorithmic}
\usepackage{stmaryrd}
\usepackage[a4paper,top=3cm,bottom=2cm,left=3cm,right=3cm,marginparwidth=1.75cm]{geometry}

\usepackage{amsmath}
\usepackage{graphicx}
\usepackage{amsfonts}
\usepackage{amsthm}
\usepackage{amssymb}
\usepackage{hyperref}
\usepackage{paralist}
\usepackage{natbib}
\usepackage{sidecap}
\usepackage{mathtools}
\usepackage{authblk}

\hypersetup{
    colorlinks=true,
    linkcolor=blue,      
    citecolor=darkgreen,  
    urlcolor=blue        
}
\definecolor{darkgreen}{rgb}{0.0, .5, 0}

\newtheorem{theorem}{Theorem}[section]
\newtheorem{lemma}[theorem]{Lemma}
\newtheorem{proposition}[theorem]{Proposition}
\newtheorem{corollary}[theorem]{Corollary}

\theoremstyle{definition}
\newtheorem{definition}{Definition}[section]
\theoremstyle{remark}
\newtheorem{remark}[theorem]{Remark}

\newtheorem{example}{Example}[section]

\newif\ifcomm
\commtrue

\newif\iflong
\longtrue

\newcounter{assumption}
\renewcommand{\theassumption}{A\arabic{assumption}}

\newcommand{\beq}{\begin{equation}}
\newcommand{\eeq}{\end{equation}}

\ifcomm
   \newcommand\comm[1]{\textcolor{blue}{ #1}}
   \newcommand{\mtodo}[2]{\todo{{\bf #1}: #2}} 
   \def\here#1{{\bf $\langle\langle$#1$\rangle\rangle$}}

\else
   \newcommand\comm[1]{}
   \newcommand{\mtodo}[2]{}
   \def\here#1{}
\fi

\newcommand{\z}{\mathbf{z}}
\newcommand{\y}{\mathbf{y}}

\def\be{\begin{equation}}
\def\ee{\end{equation}}

\def\a{\mathbf{a}}
\def\b{\mathbf{b}}

\def\p{\mathbf{p}}

\def\y{\mathbf{y}}

\def\x{\mathbf{x}}

\def\1{\mathbf{1}}

\def\b{\mathbf{b}}

\title{On the admissibility of bounds on the mean \\of discrete, scalar probability distributions\\ from an i.i.d.~sample}
\author{
    Erik Learned-Miller \hspace{.5in} George Bissias\\
    College of Information and Computer Sciences\\
    University of Massachusetts Amherst
}
\begin{document}

\maketitle
\begin{abstract}
We address the problem of producing a lower confidence bound for the mean of a discrete probability distribution, with known support over a finite set of real numbers, from an i.i.d. sample of that  distribution. Equivalently, the sample can be represented by its count vector at the known support points, which follow a multinomial distribution, and the mean of the original distribution is a known linear combination of the multinomial cell probabilities. Corresponding to each total ordering of the multinomial sample space is a bound that is conditionally optimal among bounds respecting that order. We show that every admissible bound is of this form and establish a necessary-and-sufficient condition for global admissibility: the bound must be unchanged under every total order compatible with its tied values. Thus injective, conditionally-optimal bounds are admissible, while a breakable tie proves inadmissibility. For an $N$-element sample space, there are at most $(N-1)!$ admissible bounds. Furthermore, for non-trivial sample sizes and confidence levels, at least two admissible bounds exist; hence no uniformly best valid lower bound exists.
\color{black}
\end{abstract}
\newpage
\tableofcontents
\newpage
\section{Introduction}
Let $F$ be a categorical distribution over a finite set $S$ of real-valued outcomes. For example, for $S$=$\{1,2,...,6\}$, $F$ might represent the probability of each outcome of a six-sided die whose sides are numbered one through  six. A sample of $n$ i.i.d. draws from $F$ implicitly defines a multinomial distribution $F_n$ over the sample space of counts for the categories in $S$. 

In this work, we consider the problem of establishing a lower bound on the mean of such a categorical distribution defined over a subset of real numbers from a multinomial sample. We consider the specific setting in which the support $S$ of the distribution is known, but the probabilities of each categorical outcome are unknown. Equivalently, we could bound the mean of the corresponding multinomial distribution, which is simply $n$ times the mean of the categorical distribution. We choose to formulate our problem in terms of the categorical mean, due to its closer connection to traditional confidence intervals, which focus on the mean of the underlying distribution, not the mean of the multi-sample distribution.

We rigorously define the setting below, but for the purposes of the introduction, our goal for a confidence level $1-\alpha$, support set $S$, and sample size $n$ is to produce a lower bound on the mean such that the probability of error is not greater than $\alpha$ for any distribution defined over that support. In the literature on confidence intervals, such a bound is often referred to as {\em conservative}. We prefer the term {\em valid}, since in our view, a $1-\alpha$ confidence bound must have an error rate no greater than $\alpha$ to satisfy the definition. 


We revisit an idea due to \cite{Buehler1957}. The Buehler construction begins by placing all possible samples in a fixed order, and restricts the set of bound functions considered to be consistent with that order. That is, for an ordering $T$ of a sample space and samples $\x$ and $\y$ from that sample space, $B(\x)< B(\y) \iff \x <_T \y$. 
For the lower bounds considered here, the optimization for a given sample minimizes the mean over distributions subject to a constraint on the probability of the collection consisting of that sample and all those that follow it in the order. Thus, for any two samples, the collection for one is nested within the other. This simple structure facilitates analyzing validity and comparing the strength of bounds that follow the same order. We refer to the minimization above as the {\em central optimization problem}. Applying it to each sample in one of the $N!$ possible orderings of a sample space of size $N$ yields the optimal bound for that order.\color{black} Bounds of this nature are sometimes referred to as \emph{Buehler bounds} and are said to be \emph{Buehler-optimal} with respect to the order that they induce. \cite{GuerreroDavid1985Buehler} introduced the first formalization of these bounds, deriving several extensions, and \cite{LloydKabaila2003} provided a modern treatment of Buehler bounds in a more general context.


Beyond validity lie the questions of  admissibility and optimality of lower bounds. An optimal bound is one that is valid and dominates   every other valid bound (is stronger than or equal for each sample). Admissible bounds are valid bounds that cannot be uniformly dominated. Here, the literature is less complete. \cite{KabailaLloyd2000} established that when $\alpha \in (0, 1/2)$, and under certain conditions imposed on the parameter of interest, there can be no optimal bound. In the more limited setting where the parameter of interest is the difference of two binomial proportions, \cite{Wang2010} showed that a bound is admissible among a circumscribed  class of distributions when it is injective. One of our objectives is to address questions of admissibility and dominance for lower bounds on the mean.

Another objective is to present Buehler’s approach in a form accessible to practicing statisticians and engineers by rederiving classical results in an intuitive framework specialized to mean bounds for distributions with known finite support (multinomials).
Buehler's analysis is invoked explicitly in the literature, e.g., \cite{reiser1991comment}, \cite{jobe1992buehler}, \cite{willits1997series}, \cite{lloyd2007exact}, and \cite{bancal2022simple}. The ideas also commonly arise independently. 
For example, although \cite{StarkTrinomial} do not frame their method in these terms, it closely follows Buehler's construction. They first round each original observation upward to one of three values, producing a three-category sample space. They then rank the samples in this transformed space by their sample means and maximize the population mean over all three-category distributions satisfying the same order-based probability constraints used in Buehler's upper-bound construction. Because there are only three categories, this calculation remains computationally tractable. \cite{Fienberg77} also use the same fundamental optimization as the Buehler construction. They constrain the probability of a selected collection of possible sample outcomes and optimize the mean over the set of distributions satisfying that constraint. In their construction, however, the collection of outcomes is chosen separately for each observed sample rather than being derived from a fixed order of the sample space. As a result, the collections of outcomes chosen for different observations are not necessarily nested, and so Buehler's fixed-order validity and optimality results do not directly apply to their approach.
\color{black}

Using our framework, we derive the following novel results concerning lower bounds on the mean.

\begin{itemize}
    \item \textbf{Complete characterization of global admissibility.} Building on Buehler's fixed-order construction, we provide in Theorem~\ref{thm:admissibility_characterization} necessary and sufficient conditions for the global admissibility of exact lower confidence bounds for the mean of a distribution with known finite support. Specifically, we show that a Buehler-optimal bound is globally admissible exactly when it is invariant under every total order compatible with its tied values. Together with the classical result (rederived in Theorem~\ref{thm:conditional_optimality}), which shows a Buehler bound dominates any other bound inducing the same order, our results amount to a complete characterization of the entire admissible frontier under samplewise dominance.
    \item \textbf{Nonexistence of a uniformly best bound throughout the nontrivial parameter range.} For every support containing at least two values, each sample size $n \geq 2$, and all $\alpha \in (0, 1)$, we show in Theorem~\ref{thm:no_opt} that there exist at least two admissible bounds, proving that no lower confidence bound valid in this parameter regime uniformly dominates all others.
    \item \textbf{Practical sufficient conditions for admissibility and inadmissibility.} We show in Lemma~\ref{lem:injective_optimal_bounds} that an injective, Buehler-optimal bound is necessarily globally admissible, whereas Lemma~\ref{lem:breakable} shows that a breakable tie among bound values for different samples certifies inadmissibility. We also identify, in Lemma~\ref{lem:suff_cond_no_tie}, geometric conditions under which two samples must be assigned different values by any Buehler-optimal bound that places them successively. 
    \item \textbf{Finite structural reduction and combinatorial bounds.} In Lemma~\ref{lem:admissible_bounds} we reduce the global admissibility problem to a choice among a finite collection of Buehler-optimal bounds, which, according to Lemma~\ref{lem:degen_are_dom}, excludes those bounds arising from orders that do not begin with the sample comprised exclusively of the minimal support element. We further show in Lemma~\ref{cor:max_admis} that there are at most $(N-1)!$ admissible bounds on an $N$-element sample space.
\end{itemize}
\color{black}


\subsection{Some basic definitions}
As stated above, we will consider lower bounding the mean from an i.i.d. sample of size $n$ of an unknown distribution, given the support of the distribution. Without loss of generality, we shall assume that the least value of the support is $0$, and hence shall restrict our focus to distributions with support on the non-negative reals.  While we restrict our attention to distributions with finite support over a finite set of nonnegative reals, our results can be used to approximate bounds on continuous distributions by discretizing them, as long as they have support on a finite interval.

\begin{definition}[sample space]
Given a finite support set $S=\{s_1,s_2,...,s_m\}$ and a sample size $n$, a {\bf sample space} $\Omega(S,n)$ is the set of all multinomial samples of size $n$ over the support set. For the purpose of estimating the mean from i.i.d. samples, the order of elements in the sample is irrelevant, and we represent all instantiations of a given sample by the one in which the components are sorted from least to greatest. We write $\Omega(S,n)=\{\x=(x_1,x_2,...,x_n): x_1\leq x_2 \leq ... \leq x_n; \; x_i\in S \;\forall i\}.$  We also refer to the full set of samples in a sample space as a {\bf discrete simplex lattice}. 
\end{definition}

\begin{example}
    The support set for a typical six-sided die 
    would be $S=\{1,2,...,6\}$. Suppose one rolls such a die 5 times and obtains the values $3,4,2,6,2$. We represent this as the sorted multinomial sample $\x=(2,2,3,4,6).$
\end{example}

\begin{example}
\label{ex:omega}
    Consider a support set $S=\{0,1,3\}$ and a sample size of $n=4$. Then the induced sample space $\Omega(S,n)$ is 
    \begin{eqnarray*} 
    &&\{(0,0,0,0),(0,0,0,1),(0,0,0,3),(0,0,1,1),(0,0,1,3),(0,0,3,3),(0,1,1,1), (0,1,1,3),\\
    &&(0,1,3,3),(0,3,3,3),(1,1,1,1),(1,1,1,3),(1,1,3,3),(1,3,3,3),(3,3,3,3)\}.
    \end{eqnarray*}
\end{example}
Below, we visualize the sample space as a simplex lattice in two dimensions.
\begin{center}
\includegraphics[width = 0.5\textwidth]{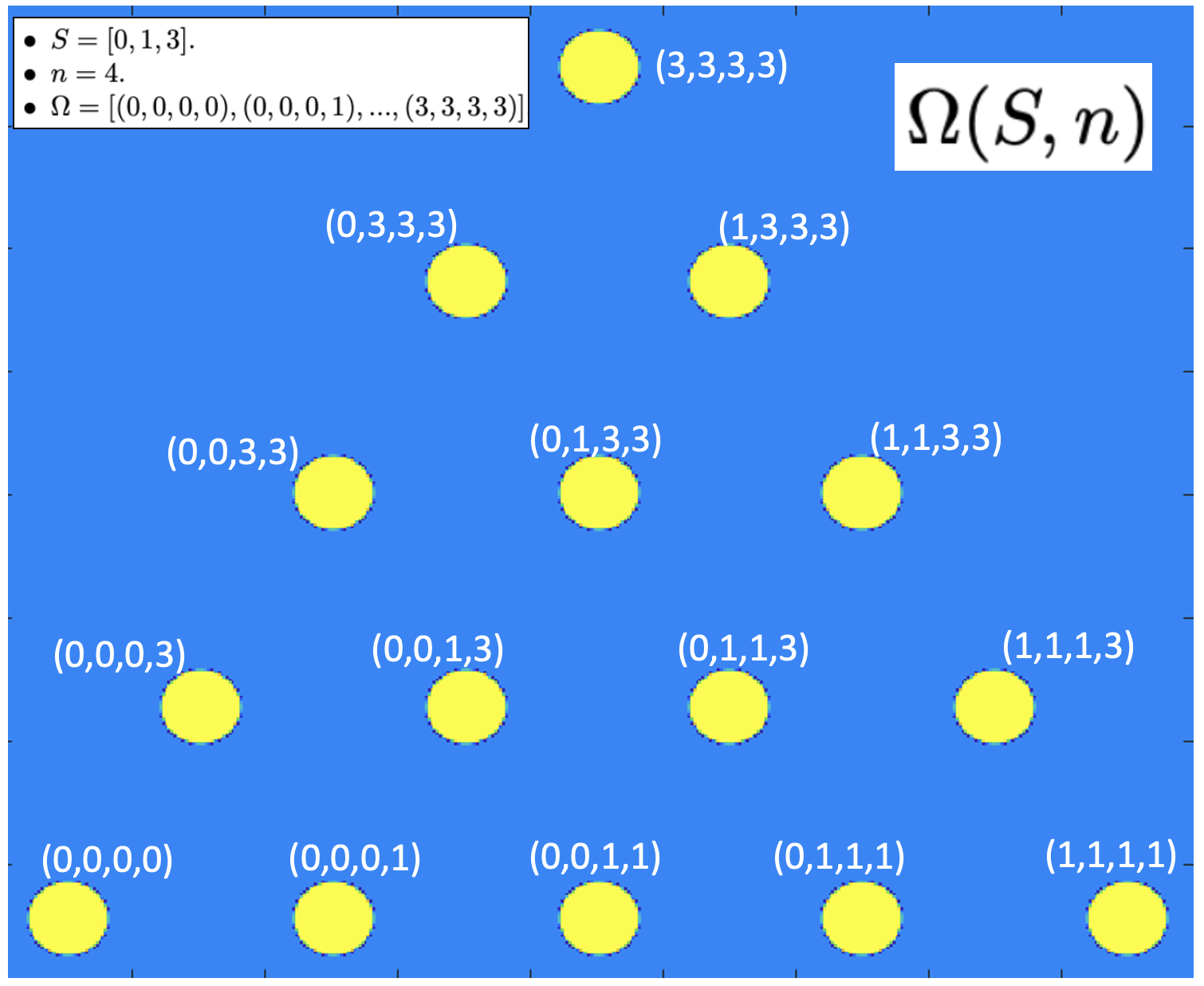}
\end{center}
\begin{definition}[multinomial likelihood function]
Consider a sample $\x$ from a multinomial distribution. Let $\hat{\x}$ be a count vector, with one entry for each possible outcome. For example, for a multinomial distribution over the values $[0,3,5,8]$ with a sample size of $7$, the vector $\hat{\x}=(4,2,0,1)$ indicates there are four 0's, two 3's, zero 5's, and one 8. It would correspond to the vector $\x=(0,0,0,0,3,3,8)$.

We can represent the probability of this outcome as a function of the parameter vector $\p=(p_1,p_2,p_3,p_4)$ of the underlying categorical distribution.
$$
L(\p|\x) \equiv \Pr(\x|\p)=\binom{n}{\hat{x}_1,...,\hat{x}_4} \prod_{i=1}^4 p_i^{\hat{x}_i},$$
where the factor in parentheses is the multinomial coefficient.
\end{definition}

\begin{example}
    For a support $S$ of size $3$, we can visualize the multinomial likelihood function for a particular sample as a function over the simplex representing the set of probability distributions over $S$. For example, for $\x=(0,0,1,3)$, we have 
    $$
    L(\p|\x) =  \binom{4}{2,1,1}p_1^2\;p_2^1\; p_3^1. 
    $$
    The left side of Figure~\ref{fig:multinomial_likelihood} shows this function over the simplex in two  dimensions.

\begin{figure}[ht]
\begin{center}
\includegraphics[width = 0.45\textwidth]{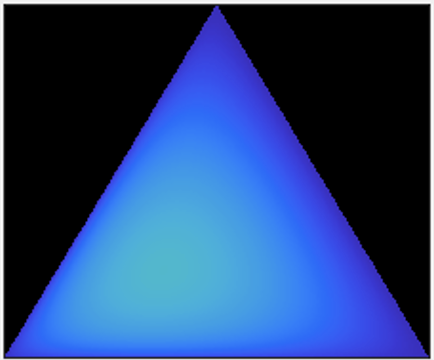}
\includegraphics[height = .351\textwidth]{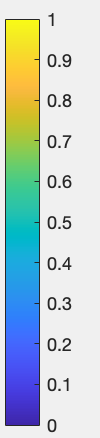}
\includegraphics[width= 0.45\textwidth]{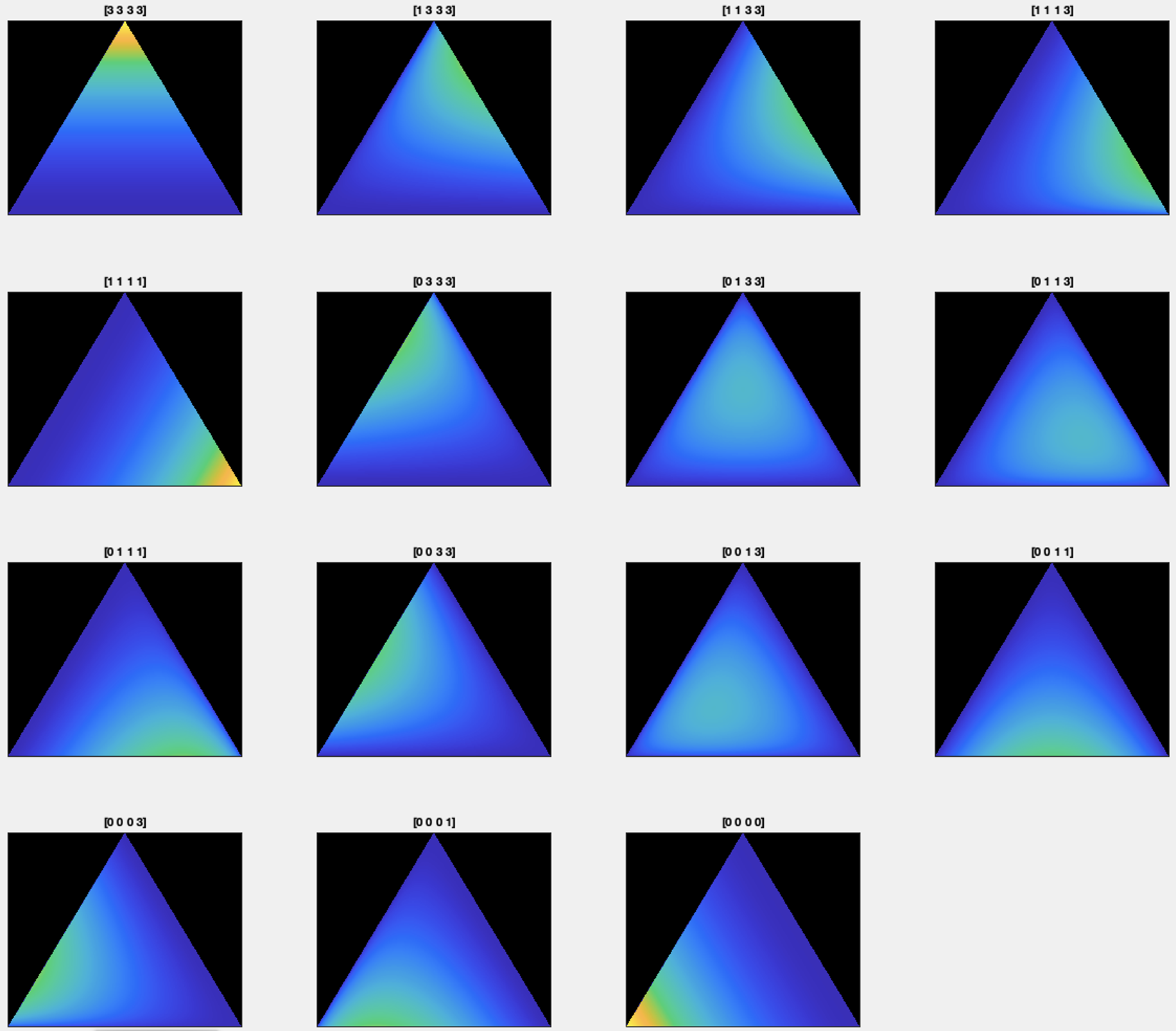}
		\caption{{\bf Left.} The multinomial likelihood function $L(\p|\x)=\Pr(\x|\p)$ for $\x=(0,0,1,3)$. Each point in the simplex gives the probability of obtaining that sample as an i.i.d. sample of size 4 from the corresponding probability distribution. The scale of the probabilities is shown on the right. Starting from the top of the triangle, and going clockwise, the three distributions at the corners of the triangle represent the distributions with all of their mass on the outcomes of 3, 1, and 0 respectively. {\bf Right}. The set of likelihood functions for each of the 15 samples in the discrete simplex for the sample space $\Omega$ of Example~\ref{ex:omega}.
        }
  		\label{fig:multinomial_likelihood}
		 \end{center}
\end{figure}
\end{example}

\begin{figure}
\begin{center}
\includegraphics[width = 0.3\textwidth]{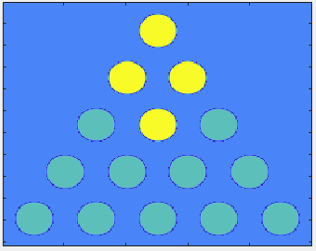}
\includegraphics[width = 0.3\textwidth]{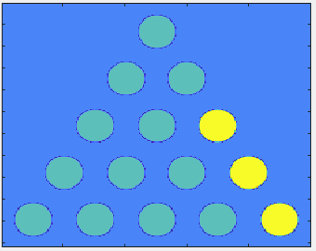}
\includegraphics[width = 0.3\textwidth]{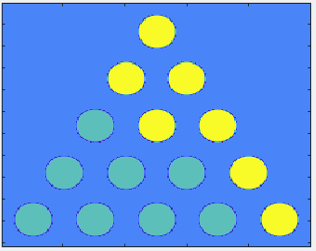}
\includegraphics[width = 0.04\textwidth]{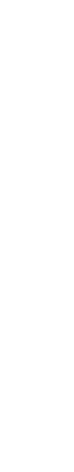}
\includegraphics[width = 0.30\textwidth]{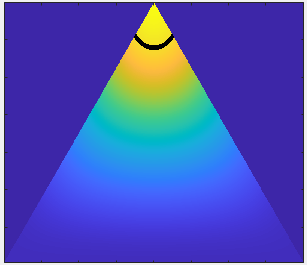}
\includegraphics[width = 0.30\textwidth]{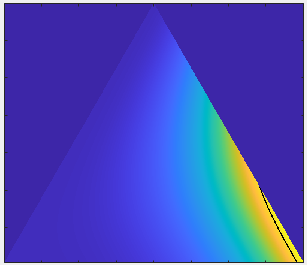}
\includegraphics[width = 0.30\textwidth]{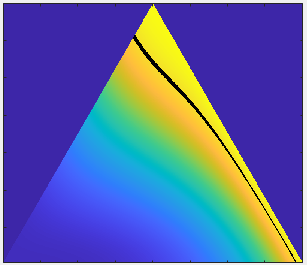}
\includegraphics[width = 0.04\textwidth,height=.26\textwidth]{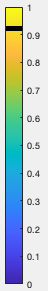}
\end{center}
		\caption{{\bf Top.} Various subsets, indicated by the yellow circles, of the full sample space. {\bf Bottom.} The multinomial likelihoods of the subsets shown at the top. The black lines illustrate a particular isocontour of the probability function, which is relevant to the central optimization problem discussed below. Notice that each subset likelihood is a sum of some subset of the sample likelihood functions shown on the right side of Figure~\ref{fig:multinomial_likelihood}.
        }
  		\label{fig:multinomial_set_likelihood}
\end{figure}

\begin{definition}[multinomial likelihood function for sets]
In addition to defining the likelihood function for a specific multinomial sample, we can extend this definition to any subset $\Omega_s$ of a sample space $\Omega$. We simply define the likelihood function for a subset $\Omega_s$ as
\begin{eqnarray}
\label{eq:multinomial_set}
L(\p|\Omega_s)=\sum_{\x\in \Omega_s} L(\p|\x).
\end{eqnarray}
This represents the probability of obtaining an outcome that is in $\Omega_s$ for each distribution in the simplex. Figure~\ref{fig:multinomial_set_likelihood} shows three examples of multinomial likelihood functions for subsets of a sample space.
\end{definition}

\begin{definition}[lower bound with specified support]
Let $S$ be a discrete, finite support set, $ n >0$ be a positive integer sample size, and $\Omega(S,n)$ be the associated sample space.
A \emph{lower bound with specified support} $B(\x; S)$ is a map $B: \Omega(S,n) \rightarrow \bar{\mathbb{R}}$, where $\bar{\mathbb{R}} \coloneqq \{-\infty\} \cup \mathbb{R} \cup \{\infty\}$. 
\end{definition}

Support set $S$ and sample size $n$ determine the sample space where the bound is defined and are therefore fixed when discussing a particular bound. Similarly, when considering \emph{validity} at confidence level $1-\alpha$, the confidence level will be treated as
fixed rather than as an argument of the bound function.

\begin{definition}[bound correctness for a particular sample and distribution]
Let $F$ be a distribution with mean $\mu$ over a support set $S$ and let $B$ be a lower bound on the sample space $\Omega(S, n)$. Bound $B$ is deemed {\bf correct} for a  sample $\x$ and the distribution $F$ if $\mu \geq B(\x)$. Otherwise the bound is {\bf incorrect} or {\bf erroneous} for $\x$. 
\end{definition}

Notice that correctness for a particular sample and distribution is determined solely by comparing $B(\x)$ with $\mu(F)$. In particular, the notion of correctness does not depend on a confidence level $1-\alpha$. We emphasize that we use the terms {\bf correct} and {\bf incorrect} (or {\bf erroneous}) to refer to the behavior of a bound for a particular sample $\x$ \textit{with respect} to a single distribution $F$, not its behavior on a sample space or an entire distribution. We reserve the terms {\bf valid} (or {\bf conservative}) and {\bf invalid} (or {\bf non-conservative}) for the behavior of a bound with respect to an entire sample space. In particular, we do not use the terms {\em valid} or {\em invalid} to describe the result of a bound on a single sample $\x$. 

\subsection{Error sets and validity}
For a given multinomial distribution and bound function, 
the error set of the bound is simply the subset of the multinomial sample space for which the bound is incorrect. One of the goals of this work is to characterize the entire space of possible bounds, and various subsets of those bounds, such as valid bounds, admissible bounds, and so on. We begin by characterizing the possible error sets that a bound can produce for different sample spaces. 

\begin{definition}[error set and valid set of a bound for a distribution and sample size]
\label{def:error_set_of_valid_bound}
Given a distribution $F$ with mean $\mu(F)$ over a support set $S$ and a sample size $n$, let $\Omega(S,n)$ be the sample space. Let $1-\alpha$ be the confidence level. Then the {\bf error set} of a lower bound $B$ for the distribution and the sample size is the set $E$ defined as 
$$
E(F,n) = \{\x\in \Omega(S,n): B(\x) > \mu(F)\},
$$  
that is, the subset of the sample space for which the bound is erroneous. The complement of this set, i.e., the subset for which the bound is correct, is termed the {\bf valid set}.
\end{definition}

\begin{definition}[validity of a bound for a distribution, sample size, and confidence level]
Let $F$ be a distribution over a support set $S$. Let $n$ be a sample size and $\Omega(S,n)$ be the induced sample space. Let $1-\alpha$ be a confidence level. Let $E(F,n)$ be the error set for the bound. We say that the bound $B$ is valid for the distribution $F$, the sample size $n$, and the confidence level $1-\alpha$ if 
$$
\Pr\nolimits_F(E) = \sum_{\x\in E} \Pr\nolimits_F(\x) \leq \alpha.
$$
\end{definition}

\subsubsection{Probability simplices and validity over a support set}
In the setting where the distribution's support is specified,\footnote{For  an example of an unknown distribution with known support, consider an unfair six-sided die. The outcomes are directly observable, but their probabilities (and hence the mean) are not.} we know that the generating distribution belongs to a particular set of distributions characterized by the probability simplex over this support set. Thus, the task of evaluating a bound over a support set is to analyze its performance with respect to each distribution in such a simplex. We now provide some more precise definitions.  
\begin{definition}[open probability simplex and closed probability simplex]
Let $S$ be a finite support set. Let $\mathcal{G}(S)$ represent the open set of probability distributions that assign strictly positive probabilities to each element of the support. We refer to this as the \textbf{open probability simplex} or simply the \textbf{open simplex}. Adding to this set all of the distributions that have support that is a subset of $S$, representing the boundary of the open simplex, we obtain the closed set $\mathcal{F}(S)$, the \textbf{closed probability simplex}, or \textbf{closed simplex}. 
\end{definition}
In some contexts, we will refer to a set of distributions whose support is strictly \textbf{equal} to a support set $S$ (open simplex). In others, we will refer to sets of distributions whose support is any \textbf{subset} of the support set (closed simplex). We will specify this when it is relevant and not clear from context.

\begin{definition}[validity of a bound for a support set]
Given a support set $S$, a sample size $n$ and a confidence level $1-\alpha$, we say that a bound is valid with respect to the support set, sample size, and confidence level if it is valid for each distribution $F\in \mathcal{F}(S),$ the closed probability simplex over $S$.
\end{definition}

\begin{definition}[validity of a bound for a set of distributions]
\label{def:validity_for_dists}
For a set of probability distributions $\mathcal{F}$, a confidence level $1-\alpha$ and a sample size $n$, a bound is called {\bf valid}, or alternatively, {\bf conservative}, if  it is valid for each  $F\in\mathcal{F}$. 
\end{definition}

\subsubsection{The structure of bound error sets over the simplex}
In this section, we pose the following question. Consider a specific bound function $B(\x)$. For the set of distributions defined over a particular sample space, how many different error sets are possible? Let $S$ be a discrete support set with least element 0. Let $\Omega(S,n)$ be a sample space for samples of size $n$.  Let $N$ be the number of elements in $\Omega$.
Now consider a bound $B(\x)$ that maps each $\x \in \Omega$ to a nonnegative real. 

\begin{definition}[injective and many-to-one bounds]
We say that a bound function $B(\x)$ is \textbf{injective} or {\bf 1-to-1} with respect to a sample space $\Omega$ if $B(\x)$ takes a distinct value for each $\x \in \Omega$. If it is not injective, then it is {\bf many-to-one} or {\bf non-injective.}
\end{definition}

Next we consider various {\em orderings} of the elements of a sample space. In order to describe such orderings, it is useful to have a default ordering and to define other orderings as permutations of this default ordering. For the default ordering, we use the classical {\em lexicographic} ordering of samples. If $\x$ and $\y$ are two vectors, then 
$\x<\y$ according to a lexicographic ordering if and only if, for the component with least index in which they disagree (call it the $i$th component), $x_i<y_i$.

\begin{definition}[sample ordering]
For a sample space $\Omega$ with $N$ elements, a \textbf{sample ordering} or \textbf{sample order} $T=(t_1, t_2,..., t_N)$ specifies an ordering $(\x_{t_1},\x_{t_2},...,\x_{t_N})$ of the samples in $\Omega$. Here, each $t_i$ is an index of the default lexicographic ordering, with the stipulation that no two components of $T$ are equal. That is, $T$ is a permutation of the lexicographic ordering.
\end{definition}

\begin{definition}[order-consistent bound]
\label{def:order_consistency}
Let $T$ be a sample ordering for a sample space $\Omega$. A bound $B$ on $\Omega$ is {\bf  consistent} with the order $T$ (or \textbf{order-consistent}) if and only if it satisfies
$$
B(\x_{t_1})\leq B(\x_{t_2}) \leq ... \leq B(\x_{t_N}).
$$
\end{definition}
\begin{remark}
\label{remark:injective}
    An injective bound is consistent with a unique sample order: the order of the bounds of the samples.
\end{remark}

\begin{remark}
\label{remark:noninjective}
    A non-injective bound is consistent with at least two sample orders.  
    For $k \geq 1$, consider a set $\{\x_{t_i},\x_{t_{i+1}},...,\x_{t_{i+k-1}}\}$ of samples with the same bound value under a non-injective bound. There are $k!$ ways of ordering these samples, so there must be at least $k!$ sample orderings compatible with such a bound. In general, for a non-injective bound with $c$ clusters of tied samples, whose cluster sizes are given by $k_1,k_2,...,k_c$, the number of orderings consistent with the bound is $\prod_{i=1}^c k_i !$. 
    As an example, consider a bound that maps the samples $\x_a$, $\x_b$, and $\x_c$ to the numbers $0$, $1$, and $1$. This mapping is consistent with two sample orders:  $\x_a\leq \x_b \leq \x_c$ and $\x_a\leq \x_c \leq \x_b$. 
\end{remark}

We now consider error sets for injective bounds.

\begin{lemma}
\label{lem:numberOfErrorSets}
Let $\Omega$ be a sample space with $N$ elements and $1-\alpha$ a confidence level. 
Let $B$ be an injective bound over $\Omega$. Then for all distributions $\mathcal{F}$ over $\Omega$, there are at most $N+1$ possible error sets for $B$. 
\end{lemma}
\begin{proof}
Let $b_1< b_2< ...<b_N$ be the values of the bounds over the sample space, sorted from least to greatest. Because the bound is injective by assumption, the inequalities are strict.

Let $\mathcal{F}$ be the closed simplex over $S$. For the distributions in $\mathcal{F}$, the only `feature' of a distribution that affects the correctness of the bound $B$ is the mean itself.\footnote{In this analysis, we adopt the convention that a bound may be considered ``correct'' or     ``incorrect'' for a distribution $F$ even if the sample has probability 0 under that distribution, which may occur for distributions that are on the boundary of the simplex.} Each possible mean value divides the sequence of bounds $b_1, \ldots, b_N$ into two subsequences: those less than or equal to the mean (correct), and those greater than the mean (errors). For a sequence of length $N$ with $N$ unique values, there are only $N+1$ positions at which it can be divided by comparison to an arbitrary real number. This of course leads to a decomposition of the sample space $\Omega$ into an error set $E$ and a valid set $V$ such that $E \cup V = \Omega$. Thus, for a given bound and sample space, there are only $N+1$ possible error sets, consisting of $0$ errors, $1$ error, $...$, up to $N$ errors. 
\end{proof}

Next consider non-injective bounds. 
\begin{lemma}
Let $\Omega$ be a sample space with $N$ elements and $1-\alpha$ a confidence level. 
Let $B$ be a non-injective bound over $\Omega$. Then for all distributions $\mathcal{F}$ over $\Omega$, there are fewer than $N+1$ possible error sets.
\end{lemma}
\begin{proof}
    In this case, some of the bound values are equal, giving 
$$
b_1\leq b_2 \leq ... \leq b_N.$$

For any pair of bound values that are equal, they must either both be correct or both be incorrect, thus reducing the number of possible splits into error sets and valid sets. Thus the number of error sets is strictly less than $N+1$ for non-injective bounds. In particular it is at most $U+1$, where $U$ is the number of unique values in the range of the bound.
\end{proof}

\subsection{Methods for comparing bounding functions}
Let $\Omega$ be a sample space over a support set $S$ and a sample size $n$.  Let $A$ and $B$ be two bounding functions which produce lower bounds for each sample $\x$ in the sample space. Of course, among multiple bounds, higher lower bounds are stronger for a given sample. But bounds are maps from a sample space to a set of bounds, one for each sample. How do we compare such bound functions in their entirety?

 We define three distinct ways of putting an order on the \textit{valid} bounding functions over a sample space. We refer to these as {\bf sample-aligned} comparison, {\bf rank-ordered} comparison, and \textbf{expected value} comparison. The first two methods lead to partial orders on valid bound functions, while the third maps each bound on a sample space to a scalar, thus leading to a total preorder. (It will be a preorder rather than a total order since ties are possible.)

 In this document, we shall restrict almost all of our analyses to the sample-aligned comparison, but we feel it is important to acknowledge the variety of ways such bound functions may be compared. In particular, many of our results presented below are particular to the sample-aligned comparison metric, and do not necessarily hold for other metrics.

\subsubsection{Sample-aligned comparison}
Let $\x_1,\x_2,...,\x_N$ be the samples from a sample space $\Omega$, where the samples are in some order specified by a total order $T$. We define a {\em sample-aligned partial order} on bounding functions as follows. For lower bound functions $A$ and $B$, we say that
$$
A\leq_{SA} B \iff A(\x_i) \leq B(\x_i),\; \forall i.$$

That is, $A$ is less than or equal to $B$ if and only if bound $A$ produces a smaller or equal value for every sample. 
Because this will be our default mechanism for comparing bounds, we will drop the $SA$ from the comparator $\leq_{SA}$ and simply use $\leq$ when the context is clear.

\subsubsection{Rank-order comparison}
Consider two lower bounds $A$ and $B$. Let $\b=(b_{(1)},b_{(2)}, ..., b_{(N)})$ be the sequence of bounds for a finite sample space $\Omega$ ranked in order from the least bound value to the greatest. Let $\a=(a_{(1)},a_{(2)}, ..., a_{(N)})$ be the ranked sequence of bound values for bound $A$.  We define a preorder on bounds $A$ and $B$ that is rank-ordered according to 
$$
A \leq_{RO} B \iff a_{(i)} \leq b_{(i)}, \; \forall i.
$$
That is, this comparison focuses on whether the $k$th ranked value of one bound is greater or less than the $k$th ranked value of another bound. 

\subsubsection{Expected value comparison}
Consider two valid bounds $A$ and $B$ over the same sample space. Let $F$ be a particular categorical distribution and let $X$ be a random element of the sample space distributed according to $F_n$ (the multinomial associated with $F$). Consider the expected value of  bound $A$:
$$
E_F[A(X)]=\sum_{\x\in \Omega} A(\x) \Pr\nolimits_F(\x).
$$
Generally, we might say that a valid bound $B$ (at the same confidence level) such that $E_F[B(X)]>E_F[A(X)]$ is a {\em stronger bound} with respect to the distribution $F$.  Of course, for some other distribution $G$ over the same sample space this relationship could be reversed, with $A$ stronger than $B$. 

One method for comparing two different bound functions across the full set of distributions for a given sample space would be to introduce a weighting or prior over the distributions in the simplex associated with the support set. Let $\mathcal{F}$ be the set of distributions on a closed simplex, and let $D(F)$ be a probability measure (e.g., a Dirichlet distribution) over these distributions. Then the expected value of a bound function $B$ over the entire simplex and the sample space would be given by
$$
E[B]=\int_{F\in\mathcal{F}} E_F[B(X)] D(F) \; dF.
$$

Given a particular prior $D(F)$ over distributions, one may wish to choose a valid bound that optimizes (maximizes) the expected value of the bound with respect to $D$.  Since this measurement of bound quality maps each bound to a scalar, it yields one or more optimum bounds with respect to this measure. We will return to an analysis of this comparison style in future work.

\subsection{Dominance, admissibility, and optimality}
The concepts of dominance, admissibility, and optimality are used to describe the relative strength of bound functions, either pairwise, or with regard to a larger set.  Because these concepts are defined with respect to a comparison metric like those defined above, which bounds are stronger will vary depending upon the underlying metric. For example, a bound that is admissible with respect to one metric may not be admissible with respect to another.

Thus, we will need to define these terms in the context of each comparison metric. However, since the sample-aligned metric is the most common, we will take the notions of dominance, admissibility, and optimality to refer to those under this metric when not otherwise stated.
In the next two sections, we introduce the basic form of our bound, and analyze its properties with respect to the sample-aligned metric.

\section{Bounds under the sample-aligned comparison metric}
\label{sec:sample_aligned_comp}
Our primary goal is to produce bounds that are ``as good as possible'' in some sense. Here we focus on this goal with respect to the sample-aligned comparison metric. 
To begin, we define the concepts of dominance, admissibility, and optimality with respect to this metric.

\subsection{Dominance, admissibility, and optimality under a sample-aligned comparison}
\begin{definition}[domination of one valid bound by another]
Let $\x_1,\x_2,...,\x_N$ be the elements of a sample space $\Omega$. Let $B_1$ and $B_2$ be two valid lower bounds over $\Omega$. We say that $B_1$ {\bf strictly dominates}  $B_2$, or simply {\bf dominates} for brevity, if $B_1(\x_i) \geq B_2(\x_i)$ for all $i$, with strict inequality for at least one $i$. We say that $B_1$ {\bf weakly dominates} $B_2$ if the weak inequalities hold for all $i$, without requiring a strict inequality for any $i$. We write $B_1 \geq B_2$ for weak dominance and $B_1>B_2$ for (strict) dominance.\color{black}
\end{definition}

\begin{definition}[admissibility with respect to a set of bound functions $\mathcal{B}$]
\label{def:admissibility}
Let $\mathcal{B}$ be a set of lower bound functions. We say that a specific bound $B\in \mathcal{B}$ is \textbf{admissible} with respect to $\mathcal{B}$ if it is \bf{valid} and there exists no \bf{valid} bound $C\in\mathcal{B}$ such that $C$ dominates $B$.\footnote{In other areas such as economics and the study of algorithms, the term {\bf Pareto optimality} is often used to describe the property of admissibility.}
\end{definition}

\begin{definition}[optimality of a bound]
\label{def:opt}
    Bound $B$ is \emph{optimal} with respect to  a set of bound functions $\mathcal{B}$ if it is valid and dominates every other valid bound in $\mathcal{B}$. It is immediate that an optimal bound must be the only admissible bound in $\mathcal{B}$. 
\end{definition}

\begin{remark}
     We will show later, in Proposition~\ref{prop:unique_implies_opt}, that for the class of bounds considered in this paper, if an admissible bound is unique, it must also dominate all other bounds.
\end{remark}
\color{black}

\subsection{Bounds consistent with specific sample orderings}
Previously, we remarked that when a bound is applied to all of the elements of a sample space, the resulting set of bounds are consistent with a single total order if the bound is injective (Remark~\ref{remark:injective}) and for two or more orderings if the bound is non-injective (Remark~\ref{remark:noninjective}).
We now consider the reverse.  In particular, given a sample ordering, we consider the set of all possible bound functions which are consistent with that order.

We aim to establish the following:
\begin{enumerate}
    \item Consider a set of bounds $\mathcal{B}$ over a sample space $\Omega$ with $N$ elements. Let 
    $\mathcal{T}=\{T_1,T_2,...,T_{N!}\}$ be the complete set of $N!$ sample orders over $\Omega$.  For the order-consistent bounds  $\mathcal{B}_{T_j}$ over $\Omega$, there is a unique optimal bound $B^*_{T_j}\in \mathcal{B}_{T_j}$ with respect to the subset of bounds $\mathcal{B}_{T_j}$. It is important to note that, except in certain trivial cases with very small sample spaces, $B^*_{T_j}$ {\em cannot} be an \textit{optimal} bound with respect to the larger set of bounds $\mathcal{B}$ (See Theorem~\ref{thm:no_opt}). 
    \item For a bound to be admissible with respect to the full set $\mathcal{B}$, it must be optimal for some subset $\mathcal{B}_T$. Otherwise, it is dominated by some other bound (the optimal one for that sample order) and cannot be admissible. Hence, the number of admissible bounds over a sample space can be no greater than $N!$, the number of sample orderings.
    \item As we shall see, given a sample order $T$ over a sample space $\Omega$, the optimal bound consistent with that sample order can be specified one sample at a time. That is, a specific order-consistent bound for a sample $\x$ is specified by a simple optimization problem that is {\em functionally independent of the bounds for the other samples}. This dramatically simplifies the definition and computation of these  bounds. 
\end{enumerate}
This approach of analyzing order-consistent bounds will lead to a number of results characterizing the space of possible bounds.

\subsection{How conditioning on a sample space ordering makes bound specification easy}
Consider a sample order $T=(t_1,t_2,...,t_N)$ which, for a specific sample space, specifies constraints on the order of bounds on those samples from least to greatest. As mentioned above, this sample ordering associates an index with each sample $\x$ giving the position of its bound among the sorted set of bound values for the sample space.\footnote{In practice, these indices refer to a fixed reference ordering, which we take to be the lexicographic ordering of the samples. That is, an ordering $T$ is specified by giving a permutation of the standard lexicographic ordering.}  For any bound that obeys this ordering, this allows us to specify an optimization problem whose solution is the optimal bound for each sample conditioned on the sample order. Together, these sample-specific bounds specify an optimal bound function over the sample space, with respect to other bounds that obey the same total order. We proceed as follows.

Given a total order $T$, suppose we wish to specify a lower bound value $B(\x_{t_k})$ for the $k$th element in the order. That is, the total order specifies that, whatever the bound for $\x_{t_k}$, it should be greater than or equal to the bounds for $\x_{t_1},...,\x_{t_{k-1}}$ and less than or equal to the bounds for $\x_{t_{k+1}},...,\x_{t_N}$.

If the bounds are consistent with the sample order, then we have the following immediate results. Let $\mu(F)$ represent the mean of a distribution $F$:
\begin{itemize}
    \item For a particular distribution $F$, if $B(\x_{t_k})>\mu(F)$, then for all $j \geq k, B(\x_{t_j})>\mu(F)$, and all such lower bounds are erroneous with respect to the distribution $F$.
    \item Similarly, if $B(\x_{t_k})\leq \mu(F),$ then for all $j \leq k, B(\x_{t_j})\leq \mu(F)$, and all such lower bounds are correct with respect to the distribution $F$.
    \item Note that if $B(\x_{t_{k-1}})$ is erroneous for a given distribution $F$, then $B(\x_{t_k})$ must be erroneous, and it is irrelevant how it is set (as long as it respects the sample order $T$). Thus, for the purposes of deciding which set of bounds are correct or incorrect for a particular distribution $F$, the setting of $B(\x_{t_k})$ is relevant only when it is the least erroneous bound.
\end{itemize}

\subsubsection{Error sets and upper sets}
The observations above describe for a particular bound function which samples {\em \textbf{must}} produce errors for a distribution $F$, given that certain other samples produce errors for the same distribution. We formalize these ideas with two closely related terms: upper sets and error sets.

\begin{definition}[upper set of a sample and a sample order]
\label{def:upper_set}
 For a sample space $\Omega$, and a particular total order $T$, consider the subset $\Omega_k$ of samples from $\Omega$ whose position in $T$ are greater than or equal to $k$:
$$
\Omega_k=\{\x_{t_k},\x_{t_{k+1}}, ..., \x_{t_N}\}.$$
We refer to this as the {\bf upper set} of $\x_{t_k}$ with respect to the sample order $T$. 
\end{definition}

\begin{definition}[error set of a sample, a sample order, and an order-consistent bound]
Consider a sample space $\Omega$, a sample order $T$, and a bound $B$ that is order-consistent with respect to $T$. The error set of a sample $\x$ is the complete set of samples such that if the bound for $\x$ is incorrect under a distribution $F$, the bound for each of the samples in the error set must also be incorrect under the distribution $F$. That is, it is the set of samples whose bound is greater than or equal to the bound for the given sample:
$$
E(\x) \equiv \{\y\in \Omega: B(\y)\geq B(\x)\}.$$
\end{definition}

Notice that for injective bounds the error set of a sample will be equivalent to its upper set, since the only samples with bounds greater than or equal to a given sample must be elements with equal or higher indices in the sample order.  However, for non-injective bounds, it is possible that an error set will contain samples in addition to the upper set, since samples with lower indices in the sample order may have the same bound value as a given sample.  Thus, in general, the error sets of a sample are supersets of the upper sets of the sample.

We shall develop bounds based on an optimization over the {\em upper sets} of each sample. We will first argue that such a procedure produces order-consistent bounds when the resulting optimization results are all distinct, i.e., that the resulting bound function is injective. Next we will demonstrate that such a procedure produces a valid bound even when the resulting bound function is {\em not injective}. 

\subsubsection{The central optimization problem}
Let $\mathcal{F}$ be the closed simplex of probability distributions over a particular support $S$. And let $\Omega(S,n)$ be a sample space over $S$ with sample size $n$. Let $T$ be a total order that induces a set of upper subsets $\Omega_i, 1\leq i \leq N,$ on $\Omega$. Let $\mathcal{G}_k$ be the subset of distributions in $\mathcal{F}$ such that the probability of $\Omega_k$ is greater than $\alpha$:
\begin{eqnarray}
\label{eq:gk}
\mathcal{G}_k = \{F \in \mathcal{F}: \Pr\nolimits_F(\Omega_k)  > \alpha\}.
\end{eqnarray}
In order for a bound to be valid, we must ensure that it is not higher than the mean of any distribution in $\mathcal{G}_k$. If it is, then we have, for some $F$, $\Pr_F\{\x \in \Omega: \mu(F) \ngtr B(\x)\} > \alpha$, meaning that the lower bound is invalid.

\begin{figure}[ht]
\begin{center}
\includegraphics[width=.25\textwidth]{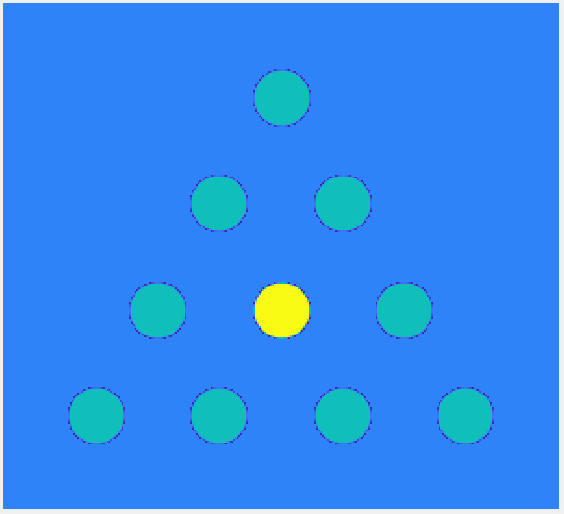}
\includegraphics[width=.25\textwidth]{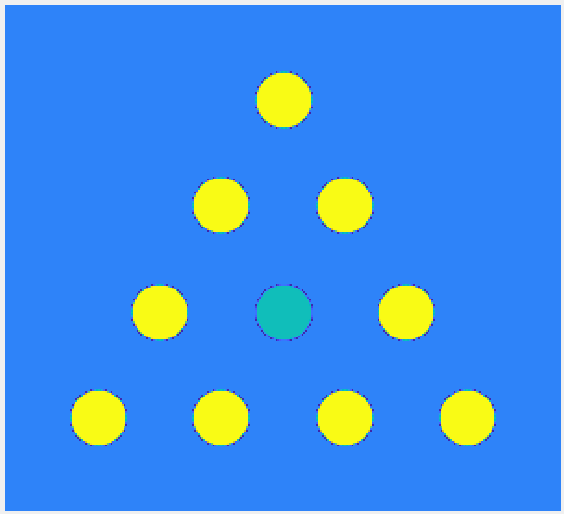}
\includegraphics[width=.25\textwidth]{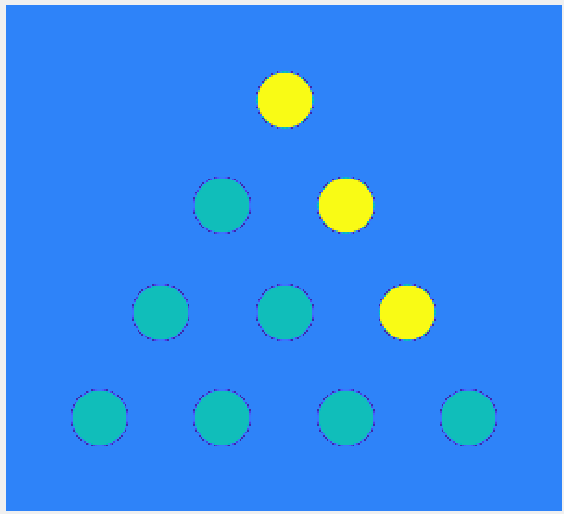}
\includegraphics[width=.03\textwidth]{figures/whitespace.png}
		\includegraphics[width = 0.25\textwidth]{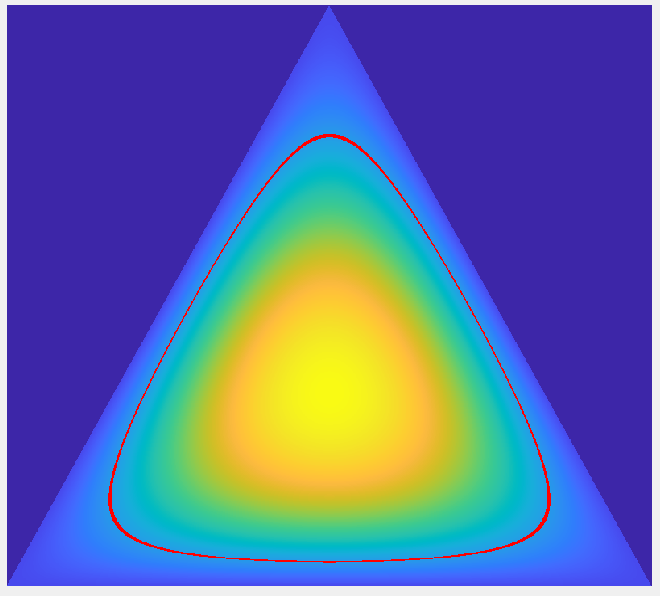}
		\includegraphics[width = 0.25\textwidth]{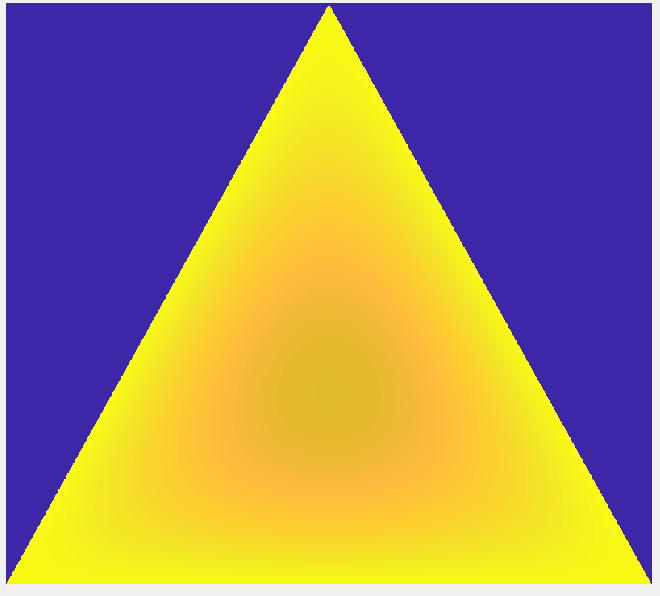}
		\includegraphics[width = 0.25\textwidth]{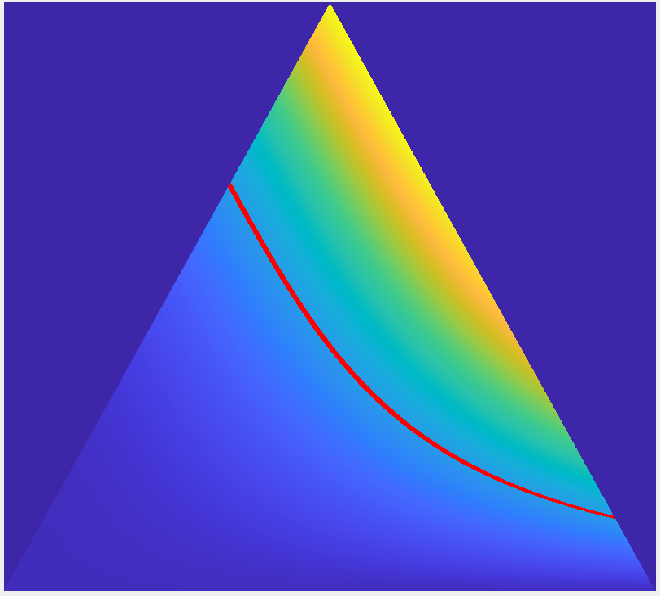}
		\includegraphics[width = 0.03\textwidth,height=1.32in]{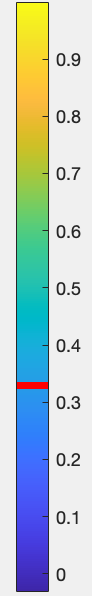}
		\caption{The top row of this figure shows three subsets (as indicated by the yellow dots) of a sample space $\Omega$ over a support set of size $3$ and a sample size of $n=3$. Below each subset, the multinomial set likelihood is shown. In each case, the red contour (if it is present) shows the the isocontour for the multinomial likelihood where the probability of the subset in the top row is approximately equal to $\alpha=0.33$. The example on the left shows how the set $\mathcal{G}_k$ of distributions with likelihoods greater than $\alpha$ forms an open set in this case. In the middle, since all distributions have probability greater than $\alpha$, the set $\mathcal{G}_k$ represents the entire simplex, and is hence closed. On the right, the set $\mathcal{G}_k$ is semi-open, having an open boundary on the bottom and closed boundaries on the top.}
        \label{fig:OpenClosedSemi}
		 \end{center}
\end{figure}

Note that the set $\mathcal{G}_k$ could be open, closed, or semi-open. To illustrate, we give examples of each case (see Figure~\ref{fig:OpenClosedSemi}).
\begin{itemize}
    \item $\mathcal{G}_k$ is open, for example, if $\Omega_k$ consists of a single sample $\x$ in the middle of the simplex lattice, and the subset of $\mathcal{F}$ with likelihood greater than $\alpha$ forms a disc-like shape around the maximum likelihood distribution. Consider the left side of Figure~\ref{fig:OpenClosedSemi}. The dark red contour shows the set of distributions $F$ such that $F(\Omega_k)=\alpha.$ The set $\mathcal{G}_k$ is  the open set inside the red contour.
    \item $\mathcal{G}_k$ is closed when $\Omega_k$ includes a large enough number of events so that the probability of being in this set is greater than $\alpha$ across the whole closed simplex. This is the case where $\mathcal{G}_k=\mathcal{F}$. See middle of Figure~\ref{fig:OpenClosedSemi}.
    \item $\mathcal{G}_k$ is semi-open when it includes part of the open simplex and part of the boundary of the closed simplex but not the entire simplex. See the right side of Figure~\ref{fig:OpenClosedSemi}. Here, $\mathcal{G}_k$ is the region above and to the right of the dark red contour (which again represents the distributions $F$ such that $F(\Omega_k)=\alpha$).
\end{itemize}

To ensure that a lower bound is no higher than the mean of any distribution in $\mathcal{G}_k$, it is sufficient to set it to the infimum of the means of the distributions in $\mathcal{G}_k$.  One drawback of this formulation is that we cannot easily refer to  a set of distributions within $\mathcal{F}$ whose means achieve the infimum. In particular, the distributions, if any, which achieve the infimum are not in $\mathcal{G}_k$ when $\mathcal{G}_k$ is open.  To simplify our analysis in future results, we form the closure of the set $\mathcal{G}_k$, and define our bound with respect to this closed set as follows:
\begin{eqnarray}
\mu^*&=& \inf_{F\in \mathcal{G}_k} \; \mu(F)\\
&=& \min_{F\in \mathcal{F}_k} \; \mu(F),
\end{eqnarray}
where 
$$
\mathcal{F}_k=cl(\mathcal{G}_k),$$
simply adding a boundary to the locally open parts of the set. We refer to $\mathcal{F}_k$ as the {\bf likely set} of distributions for the subset $\Omega_k$.

To specify a lower bound for $\x_{t_k}$, we consider two cases: where the likely set $\mathcal{F}_k$ is empty and where it is not empty.
If $\mathcal{F}_k$ is not empty we define the optimal bound to be the solution of the following optimization problem: 
\begin{eqnarray}
    \label{eq:CentralOptimizationProblem}
B^*(\x_{t_k}| \mathcal{F}_k \neq \emptyset) = \min_{F\in \mathcal{F}_k} \; \mu(F).
\end{eqnarray}
Equation~\ref{eq:CentralOptimizationProblem} is the fundamental optimization problem first studied by \cite{Buehler1957}. We refer to it here as the {\em central optimization problem} for multinomial bounds.

Notice that since the likely set $\mathcal{F}_k$ is a closed set and 
because the mean of the distributions is a continuous function over the probability simplex, the minimum will always exist. Furthermore, we can refer to the distributions in $\mathcal{F}_k$ which achieve this minimum as $\arg \min$s, and there will always be at least one such distribution.  

Algorithmically, to set the bound for the $k$th element of a sample space $\Omega$, we proceed as follows:
\begin{itemize}
    \item Form the upper set $\Omega_k$ of samples from $\Omega$ whose position are greater than or equal to $k$  in the total order.
  \item Form the polynomial representing the probability of   $\Omega_k$ (Equation~\ref{eq:multinomial_set}), as a function of the parameters of $F$. This will be a homogeneous polynomial\footnote{A {\em homogeneous polynomial} is a polynomial in which the degree of each monomial (the sum of the powers of the factors) is a constant. In this case, the degree of each monomial, such as $Cx_1^3x_2^2x_3^7$, will be equal to the sample size $n$.} with one (multinomial probability) term for each sample in $\Omega_k$. 
  \item For the likely set of distributions $\mathcal{F}_k$ (the closure of the strict superlevel set $\{F: \Pr_F(\Omega_k) > \alpha\}$), find the minimum of the means. This minimum is the lower bound on the mean. 
\end{itemize}

When $\mathcal{F}_k$ is empty, we define the lower bound to be $\infty$:
$$
B^*(\x_{t_k}| \mathcal{F}_k = \emptyset) = \infty.$$
It may seem like an odd choice to set a lower bound to $\infty$, since for every distribution over $S$ the bound for the specific sample $\x_{t_k}$ will be erroneous. But such a choice {\em can be} part of a valid lower bound over a family of distributions. 

Thus, the full specification of our order-conditioned bound, conditioned on the total order $T$, is as follows.
\begin{definition}[Buehler bound]
\label{def:cop}
    Consider a sample space $\Omega$, a sample order $T$, and sample $\x_{t_k} \in \Omega$, the $k$th lowest sample in order $T$. The \emph{exhasutive bound} on $\x_{t_k}$ with respect to $T$ is given by
    \begin{eqnarray}
    \label{eq:optBound}
    B^*(\x_{t_k}) =
    \begin{cases}
        \infty,& \text{if } \mathcal{F}_k=\emptyset \\
        \underset{F\in \mathcal{F}_k}{\min} \;\mu(F),              & \text{otherwise.}
    \end{cases}
    \end{eqnarray}
\end{definition}

Below we will prove a variety of results about this bound, but first we need a few simple lemmas to support our arguments.


\subsection{Order-conditioned bounds and conditional optimality}

In the following, for simplicity, we again assume that the least element of the support is $0$. It is simple to generalize results to other cases, but this will simplify arguments. We begin by showing that nested upper sets lead to non-decreasing bounds when those bounds arise from the Buehler bound of Definition~\ref{def:cop}. The following two results are due originally to \cite{GuerreroDavid1985Buehler}. We restate them here for continuity and also to explicitly handle bounds taking on infinite values, which was not addressed directly by their work. 

\begin{lemma}[order-consistency of Buehler bounds]
\label{lemma:incremental_sums}
Fix total order $T$. If $\x_{t_k}$ precedes $\x_{t_j}$ in order $T$, then $B^*(\x_{t_k}) \leq B^*(\x_{t_j})$. 
\end{lemma}
\begin{proof}
    By Definition~\ref{def:upper_set} we know that $\Omega_j \subseteq \Omega_k$, which implies by the monotonicity of measure that $\mathcal{F}_j \subseteq \mathcal{F}_k$. By Definition~\ref{def:cop}, the conclusion follows immediately if $\mathcal{F}_k = \emptyset$. For $\mathcal{F}_k$ non-empty, because $\mathcal{F}_j \subseteq \mathcal{F}_k$, we have
    \begin{eqnarray}
        B^*(\x_{t_k}) &=& \underset{F\in \mathcal{F}_k}{\min} \;\mu(F) \\
        &\leq& \underset{F\in \mathcal{F}_j}{\min} \;\mu(F) \\
        &=& B^*(\x_{t_j}).
    \end{eqnarray}
\end{proof}

The following definition is equivalent to that of \emph{Buehler optimality}, as defined in the literature (e.g., \cite{reiser1991comment}, \cite{jobe1992buehler}, and \cite{willits1997series}). We prefer the term \emph{conditional optimality} because it emphasizes that the optimality has limited scope.

\begin{definition}[conditional optimality of a bound]
\label{def:cond_opt}
    Consider any bound $B$ consistent with total order $T$. We say that $B$ is \emph{conditionally optimal} for $T$ if it is optimal in the sense of Definition~\ref{def:opt} with respect to every other bound consistent with $T$.
\end{definition}
\color{black}

\begin{theorem}[conditional optimality]
\label{thm:conditional_optimality}
Let $\mathcal{B}$ be the lower bound functions (both valid and invalid) over a sample space $\Omega$ and for a confidence level $1-\alpha$. For a given total order $T$, let $\mathcal{B}_T\subset \mathcal{B}$ be the order-consistent bounds (again, both valid and invalid) with respect to a total order $T$. Let $B^*\in \mathcal{B_T}$ be the Buehler bound of Definition~\ref{def:cop}. Then $B^*$ is optimal with respect to $\mathcal{B}_T$. 
\end{theorem}

\begin{proof}
The conditional optimality of $B^*$ requires optimality according to Definition~\ref{def:opt} among bounds in $\mathcal{B}_T$ as well as validity according to Definition~\ref{def:validity_for_dists}. We proceed by proving each separately.

{\bf Validity.} 
For any $F \in \mathcal{F}$, let $E^*(F, n)$ denote the error set for $F$ with respect to $B^*$, where by Definition~\ref{def:error_set_of_valid_bound}
\begin{eqnarray}
\label{eq:exhaustive_error_set}
    E^*(F,n) \coloneqq \{\x\in \Omega(S,n): B^*(\x) > \mu(F)\}.    
\end{eqnarray}
According to Definition~\ref{def:validity_for_dists} it will suffice to show that $\forall F \in \mathcal{F}$, $\Pr_F(E^*(F, n)) \leq \alpha$. 
To that end, suppose for the sake of contradiction that there exists $F \in \mathcal{F}$ such that $\Pr_F(E^*(F, n)) > \alpha$.
According to Lemma~\ref{lemma:incremental_sums}, $B^*$ itself will be in $\mathcal{B}_T$.
Therefore, from Equation~\ref{eq:exhaustive_error_set} we can see that there exists $\Omega_k$ such that $\Omega_k = E^*(F, n)$ and $B^*(\x_{t_k}) > \mu(F)$.
It follows then that $\Pr_F(\Omega_k) > \alpha$, which implies that $F \in \mathcal{F}_k$.
But since $\mathcal{F}_k$ is non-empty, Equation~\ref{eq:optBound} ensures that $B^*(\x_{t_k}) \leq \mu(G)$ for each $G \in \mathcal{F}_k$.
This is true in particular for $F$.
Thus, $B^*(\x_{t_k}) \leq \mu(F)$, a contradiction.
It must then be the case that $\Pr_F(E^*(F, n)) \leq \alpha$ for each $F \in \mathcal{F}$.
\color{black}
\medskip

{\bf Optimality for order $T$.}
Let $B \in \mathcal{B}_T$ be any valid bound. It suffices to show that 
\begin{eqnarray}
    B(\x_t) \leq B^*(\x_t)
\end{eqnarray}
for each $\x_t \in \Omega$. To that end, consider an arbitrary sample $\x_{t_k}$. If $B^*(\x_t) = \infty$, then the conclusion clearly holds. Otherwise, suppose that $B(\x_{t_k}) > B^*(\x_{t_k})$. It must then be the case that $\exists \delta > 0, B(\x_{t_k}) = B^*(\x_{t_k}) + \delta$.
And since $B^*$ is the Buehler bound of Definition~\ref{def:cop}, we must also have that $\forall \epsilon > 0, \exists G_\epsilon \in \mathcal{G}_k: B^*(x_{t_k}) + \epsilon \geq \mu(G_\epsilon)$.
Thus, by choosing $\epsilon < \delta$, we have that $B(\x_{t_k}) > \mu(G_\epsilon)$.
Of course, $\Pr_{G_\epsilon}(\Omega_k) > \alpha$ since $G_\epsilon \in \mathcal{G}_k$.
And because $B$ is consistent with $T$, it must also be the case that $B(\x) > \mu(G_\epsilon)$ for all $\x \in \Omega_k$.
This, however, would imply that the error set for $G_\epsilon$ under bound $B$ has measure greater than $\alpha$, making $B$ invalid, which is a contradiction. Thus, $B(\x_t) \leq B^*(\x_t)$ for every valid $B \in \mathcal{B}_T$ and each choice of $x_t \in \Omega$. 
\color{black}

Notice that this result does not depend upon the value of the bound for any other sample, as long as the ordering is fixed. Thus, conditioning on the order makes it possible to analyze the bound one sample at a time.
\end{proof}

\section{Injectivity, total orders, and admissibility}
\label{sec:adm_inj_total}
Now that we have established the conditional optimality of Equation~\ref{eq:optBound}, we turn our attention to the question of which conditionally optimal bounds are admissible over all orderings. We shall establish the following results.
\begin{itemize}
    \item \textbf{[Lemma~\ref{lem:admissible_bounds}].} No bound can be admissible unless it is a conditionally optimal bound.
    \item \textbf{[Lemma~\ref{lem:injective_optimal_bounds}.]} All conditionally optimal bounds that are injective are admissible.
    \item Among conditionally optimal bounds that are non-injective, there are two types: those with a \textit{breakable tie} (to be defined below), and those with \textit{unbreakable ties}. Those with unbreakable ties are admissible, while those with a breakable tie are dominated by other conditionally optimal bounds, and hence are not admissible.
\end{itemize}
We start by demonstrating that a necessary condition for a bound to be admissible is that it be one of the order-conditioned optimal bounds. 

\begin{lemma}[All admissible bounds are conditionally optimal for some total order.]
\label{lem:admissible_bounds}
Let $\mathcal{B}$ be all of the lower bounds for a discrete support set $S$, finite sample size $n$, and confidence $1-\alpha$. A lower bound $B$ cannot be admissible with respect to the set of bounds $\mathcal{B}$ unless it is conditionally optimal for some ordering.
\end{lemma}
\begin{proof}
 From Theorem~\ref{thm:conditional_optimality}, there is a unique lower bound for each total ordering over the sample space that is conditionally optimal with respect to the set of bounds that share the same total order.  If a bound consistent with total order $T$ is not conditionally optimal, it is dominated by the optimal bound for $T$, and thus cannot be admissible.
\end{proof}

\begin{corollary}
\label{cor:nfactorial}
Let $\mathcal{B}$ be the set of all lower bound functions for a discrete support set $S$, finite sample size $n$, and confidence $1-\alpha$.   Among these, there are no more than $N!$ admissible bounds, where $N$ is the number of elements in the sample space. 
\end{corollary}
\begin{proof}
There are $N!$ total orderings for a sample space with $N$ elements, and there is no more than one conditionally optimal bound per total order. Since every admissible bound must be conditionally optimal, there can be no more than $N!$ admissible bounds.
\end{proof}

\begin{proposition}
\label{prop:unique_implies_opt}
    For the finite sample spaces $\Omega$ considered in this document, let $\mathcal{B}$ be a class of bounds over $\Omega$. If $A \in \mathcal{B}$ is the unique admissible bound, then $A$ dominates every other valid bound in $\mathcal{B}$. 
\end{proposition}

\begin{proof}
    Let $A$ be the unique admissible bound and let $\mathcal{C}$ denote the set of all conditionally optimal bounds. For each $B_T \not \in \mathcal{C}$, consistent with $T$, there exists $B_T^* \in \mathcal{C}$ that is conditionally optimal for $T$, which implies by Theorem~\ref{thm:conditional_optimality} that $B_T^*$ dominates $B_T$. Thus, it suffices to show that $A$ dominates every other bound in $\mathcal{C}$.

    Because $\Omega$ is finite, allowing for finitely many total orders, there can be just one conditionally optimal bound for each order, making $\mathcal{C}$ finite as well. Now, let $B_1 \in \mathcal{C}$ be arbitrary, and build admissible bound $M$ inductively as follows. As long as $B_i$ is not admissible, choose any $B_i'$ that dominates it, and take $B_{i+1} \in 
    \mathcal{C}$ to be the conditionally optimal bound for the order induced by $B_i'$. Theorem~\ref{thm:conditional_optimality} ensures that $B_{i+1}$ weakly dominates $B_i'$ so it must strictly dominate $B_i$. Having established that $\mathcal{C}$ is finite, this induction must terminate at a finial conditionally optimal bound $M \in \mathcal{C}$. But $A$ is the only admissible bound by assumption, thus we must have that $M = A$. It follows that $A$ dominates every bound in $\mathcal{C}$.
\end{proof}
\color{black}

\subsection{Injective bounds}
Many of our results are simplified for conditionally optimal bounds that are injective with respect to a given sample space $\Omega$, i.e., bounds that do not map any two samples of $\Omega$ to the same real value. In this section, we explore properties identified with injective and non-injective bounds, and present some results specialized to one category or the other.

As discussed above, for each total order $T_j$, $(1\leq j \leq N!)$, there is a single conditionally optimal bound with respect to the subset of bounds $\mathcal{B}_{T_j}$ that are consistent with that order. We shall refer to it as $B^*_{T_j}$ and call it the conditionally optimal bound for $T_j$.  
The following lemma generalizes an earlier result from \cite{Wang2010}, which showed that, for the problem of comparing the difference of two binomial proportions, a particular recursively defined bound is admissible among a certain restricted class of bounds when that bound has no ties. 
Our lemma is also related to \cite{Wang2023}, which extends the Buehler construction to two-sided confidence intervals and proves that the limiting interval in a certain sequence is admissible when both endpoints are injective functions of the sample. This result is not directly applicable here because representing our lower bound as an interval would require setting the upper endpoint equal to $+\infty$ for every sample, making it non-injective.
\color{black}

\begin{lemma}[Injective conditionally optimal bounds are admissible]
\label{lem:injective_optimal_bounds}
 Let $B^*_T$ be a conditionally optimal bound with respect to a sample space $\Omega$ and with respect to a total order $T$. Furthermore, assume that it is injective. Then it is admissible with respect to the full family of bounds $\mathcal{B}$ over $\Omega$.
\end{lemma}
\begin{proof}
 To prove this, we must show that there is no other valid bound that dominates $B^*_T$, and it is sufficient to focus on the conditionally optimal ones. To do this, we shall show that for any injective $B^*_T$, there  exists for every other total order $U\neq T$ at least one sample $\x$ such that $B^*_T(\x)>B^*_U(\x)$, i.e., that $B^*_T$ gives the better bound. 

 Let $T=(t_1,...,t_N)$ be a total order. Let $B^*_T$ be its conditionally optimal bound, and assume that it's injective. Let $U=(u_1,...,u_N)$ be a different total order whose conditionally optimal bound $B^*_U$ is not necessarily injective.
Let $\x_{t_i}$ be the  $i$th sample according to the total order $T$, and  $\x_{u_i}$ the $i$th sample according to the order $U$. 

 Let $k$ be the least index such that $t_k\neq u_k$. Note that $k\leq N-1$ since two different lists of the same elements must differ in at least two positions. Let $\Omega^T_k=\{\x_{t_k},\x_{t_{k+1}},...,\x_{t_N}\}$ be the upper set for the sample $\x_{t_k}$ under $T$, and  $\Omega^U_k=\{\x_{u_k},\x_{u_{k+1}},...,\x_{u_N}\}$ be the upper set for the sample $\x_{u_k}$ under $U$. Because the first $k-1$ elements of each total order are equivalent, and the upper sets are the complements of the subsets of these first elements, 
 we have that $\Omega^T_k = \Omega^U_k$. Then 
 \begin{eqnarray}
 \label{eq:equalBounds}
B^*_T(\x_{t_k})=B^*_U(\x_{u_k}).
 \end{eqnarray}
 
 Let $\Omega^T_{u_k}$ be the upper set of $T$ beginning at the element $u_k$. Since $u_k$ must appear later than position $k$ in order $T$, we have that
 $$
 \Omega^T_{u_k} \subset \Omega^T_{t_k},
 $$
 so by Lemma~\ref{lemma:incremental_sums} we have that
 \begin{eqnarray*}
B^*_T(\x_{t_k}) \leq B^*_T(\x_{u_k}),
 \end{eqnarray*}
 and since by the assumption of injectivity,  $B^*_T(\x_{t_k}) \neq B^*_T(\x_{u_k})$, we conclude that
\begin{eqnarray}
\label{eq:lessThanBounds}
B^*_T(\x_{t_k}) < B^*_T(\x_{u_k}).
\end{eqnarray}
 
 Combining  Eq.~\ref{eq:equalBounds} and Eq.~\ref{eq:lessThanBounds}, we conclude that 
 $$
 B^U(\x_{u_k}) < B^T(\x_{u_k}). 
 $$
 That is, for sample $\x_{u_k}$ the conditionally optimal bound based on $T$ outperforms the conditionally optimal bound based on $U$.
 \end{proof}

\subsection{Non-injective bounds}
Lemma~\ref{lem:injective_optimal_bounds} concerns injective bounds. What about non-injective bounds? Can they also be admissible? The story is a bit more complicated here. 

To understand the landscape of non-injective bounds and which ones are admissible, we introduce the notions of breakable ties and unbreakable ties.
As mentioned in Remark~\ref{remark:noninjective}, a bound with ties is any bound over a sample space in which two or more samples are mapped to the same bound value. Of course, due to the fact that all conditionally optimal bounds must obey the same ordering as the sample ordering they are conditioned on, any ties must be from consecutive elements of the ordering. 

To gain some intuition about when bounds yield ties and when they do not, we start with a result about conditions under which the bounds for two samples adjacent in an order $T$ are guaranteed to produce different bounds. But first we need the following lemma. 

\begin{lemma}[Non-negativity of arbitrary sums of multinomial likelihoods]
\label{lem:continuity}
    Each multinomial probability is continuous, infinitely differentiable,  non-negative on the closed simplex, and strictly positive on the open simplex. Any sum of such likelihoods (that are all conditioned on the same parameters) inherits the same properties. 
\end{lemma}
\begin{proof}
    Multinomial probabilities are just polynomials with positive coefficients, so continuity and differentiability and non-negativity are trivial. Since any probability distribution in the open simplex assigns a non-zero probability to each outcome of a multinomial, every possible sample has a positive probability, meaning that all possible polynomials are positive on the open simplex. On the boundary of the simplex, the underlying distributions have at least one outcome whose probability is zero, so any multinomial comprised exclusively of samples containing at least one boundary outcome will have likelihood zero.
\end{proof}

\begin{figure}[ht]
\begin{center}
		\includegraphics[width = 0.7\textwidth]{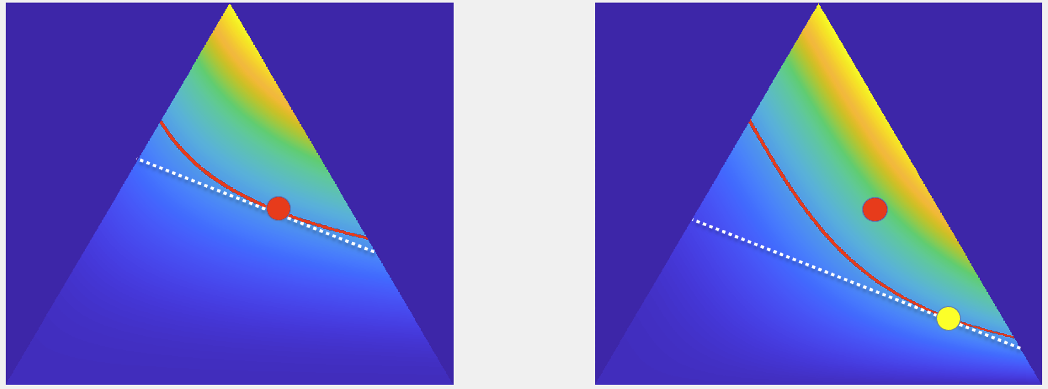}
		\caption{{\bf Left.} The left figure plots a polynomial function over the simplex corresponding to the probability of the subset $\{(1,1,3),(1,3,3),(3,3,3)\}$ of a sample space over the support $\{0,1,3\}$. The red line shows the set of distributions for which the likelihood of the set is $\alpha$ (here, $\alpha = 0.35$), and the white dotted line shows a set of distributions with the same mean (an isomean contour). The red dot shows the minimum-mean distribution in the simplex over the closure of the region where the probability of the subset is greater than $\alpha$. {\bf Right.} The right figure is similar but shows the same plot for a slightly larger subset of samples: $\{(0,1,3),(1,1,3),(1,3,3),(3,3,3)\}$. Note that the point in the left figure at which the probability was $\alpha$ now yields a substantially higher probability, and no longer represents an optimum. The optimal mean has moved to the mean of the yellow point's distribution.}
  		\label{fig:incrementalBounds}
		 \end{center}
\end{figure}

Now we return to the question of what circumstances guarantee that the central optimization problem of two different samples will not be the same. This is captured by the following lemma.

\begin{lemma}[sufficient conditions for bounds being different (non-tied)]
\label{lem:suff_cond_no_tie}
    Let $\Omega_k\neq \Omega$ be a strict subset of $\Omega$, $\x$ an arbitrary sample not in $\Omega_k$, and $\Omega_j \coloneqq \Omega_k \cup \{\x\}$. For a fixed $\alpha$, let 
    $$
     \mathcal{F}_{\min} =\underset{cl\{F \in \mathcal{F}: \Pr\nolimits_F(\Omega_k) > \alpha\}}{\arg \min} \mu(F),
    $$
    and suppose that $\mathcal{F}_{\min}$ is non-empty.
    Further suppose that there exists $F^* \in \mathcal{F}_{\min}$ that is an element of the open simplex $\mathcal{G}\subset \mathcal{F}$. That is, $F^*$ is not on the boundary of the simplex. Then we have
    $$
        B^*_j < B^*_k < \infty,
    $$
    where $B^*_k$ and $B^*_j$ are the Buehler bounds for $\Omega_k$ and $\Omega_j$, respectively.
\end{lemma}
\begin{proof} To begin, because $\mathcal{F}_{\min}$ is non-empty by assumption, $B^*_k$ is finite, and we have that for any $F \in \mathcal{F}_{\min}$, $\mu(F) = B^*_k$. 
Lemma~\ref{lemma:incremental_sums} further establishes that $B^*_j \leq B^*_k$. Thus, it suffices to rule out the possibility that $B^*_j = B^*_k$. To that end, we need only show that there exists a distribution $G$ satisfying $\mu(G) < B^*_k$ with $\Pr_G(\Omega_j) > \alpha$.

Fix any $F^* \in \mathcal{F}_{\min} \cap \mathcal{G}$, which is guaranteed to exist by assumption, and define $\mu^* \coloneqq \mu(F^*)$. 
We next show that $\Pr_{F^*}(\Omega_k) = \alpha$. Since $F^* \in \mathcal{F}_{\min}$, it must assign at least probability $\alpha$ to $\Omega_k$. Suppose for the sake of contradiction that $\Pr_{F^*}(\Omega_k) > \alpha$. Since $F^* 
\in \mathcal{G}$, there exists a neighborhood of distributions surrounding $F^*$ whose measure of $\Omega_k$ also exceeds $\alpha$. But because the mean is a non-constant linear function on the simplex, there must then exist some distribution in the neighborhood whose mean is strictly below $\mu^*$, which cannot occur given the construction of $\mathcal{F}_{\min}$. Thus, $\Pr_{F^*}(\Omega_k) = \alpha$.

Now consider adding sample $\x \not \in \Omega_k$ to $\Omega_k$ to form $\Omega_j$. Since $F^* \in \mathcal{G}$, we have by Lemma~\ref{lem:continuity} that $\Pr_{F^*}(\x) > 0$. It follows that
\begin{eqnarray*}
    \Pr\nolimits_{F^*}(\Omega_j) &=& \Pr\nolimits_{F^*}(\Omega_k) + \Pr\nolimits_{F^*}(\x) \\
    &=& \alpha + \Pr\nolimits_{F^*}(\x) \\
    &>& \alpha.
\end{eqnarray*}
But again because $F^* \in \mathcal{G}$, for some $\epsilon > 0$, there exists a distribution $F_\epsilon$, a perturbation of $F^*$, which raises the mass at one support point by $\epsilon$ and lowers the mass at a larger support point by $\epsilon$, while maintaining the property that $\Pr_{F_\epsilon}(\Omega_j) > \alpha$. Of course we also have $\mu(F_\epsilon) < B^*_k$ by construction. Therefore, according to our sufficient condition, we must have $B^*_j < B^*_k$.
\end{proof}


Let $F^* \in \mathcal{F}_{\min} \cap \mathcal{G}$ be as defined in the proof of Lemma~\ref{lem:suff_cond_no_tie} so that $\Pr_{F^*}(\Omega_k) = \alpha$ and $\mu(F^*) = B^*_k$. 
In our setting, Theorem 1 of \cite{KabailaLloyd2006} implies that there exists strict separation between the Buehler bounds for samples $\y$ and $\z$, where the designated statistic orders $\y$ before $\z$, provided that $\Pr_{F^*}(\x) > 0$ for any $\x \in \Omega$ whose Buehler bound is equal to that of $\z$. Lemma~\ref{lem:suff_cond_no_tie} provides a simple, geometric sufficient condition for their theorem to apply.
\color{black}

Figure~\ref{fig:incrementalBounds} illustrates the ideas in Lemma~\ref{lem:suff_cond_no_tie}. The left part of the figure shows the distribution which is the argmin of the central optimization problem as a red dot. Note that it is in the open simplex, not on the boundary of the simplex. The white dotted line shows an isomean contour which cannot be made any higher without intersecting the set of distributions whose probabilities are greater than $\alpha$. The right panel shows that when the central optimization problem is solved for a larger subset of the sample space, the result of the central optimization problem must move lower, since the multinomial likelihood of every point in the neighborhood of the red dot must go up. With a greater portion of the open simplex having likelihoods greater than $\alpha$ in this neighborhood, the minimum must go down.

\subsection{Bounds with ties}
\label{sec:boundsWithTies}

\begin{figure}[ht]
    \begin{center}
		\includegraphics[width = 0.99\textwidth]{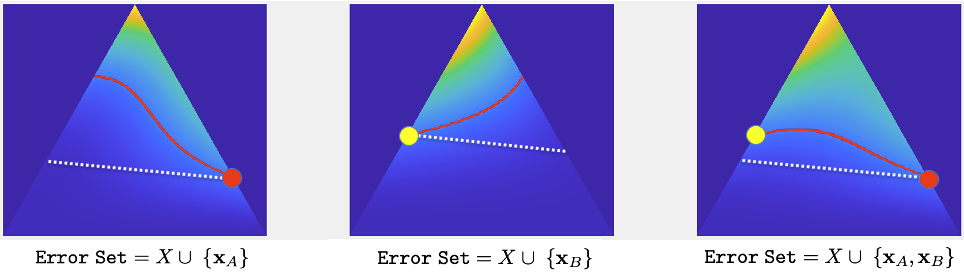}
		\caption{{\bf Left.} This figure illustrates several ideas. Each plot shows a multinomial likelihood for a different sample space subset over the same simplex (the set of probability distributions over a fixed support set). In each plot, the red curve shows the isocontour with  probability equal to $\alpha$. The white dotted lines show isocontours at the minimum mean over the closure of the region where the probability exceeds $\alpha$. The plot on the left illustrates the central optimization problem for an error set that is the union of a set of samples $X$ augmented by another single sample $\x_A$. The distribution that achieves the minimum mean is shown with the red dot. Note that this optimum distribution is on the boundary of the simplex, and hence is not in the open simplex. The rightmost figure shows the optimization problem with both $\x_A$ and $\x_B$ added to $X$. Unlike cases in which an optimum with a smaller subset occurs in the open simplex, we see that the optimal distribution is still in the same place (red dot again). Hence, for this support set and this ordering of samples ($B(\x_A) \leq B(\x_B))$, the bound has ties and hence is not injective. Note that if we had instead chosen a total order such that $B(\x_B) \leq B(\x_A)$, the bound for $\x_B$ would be strictly better, as shown by the yellow dot in the central figure, and the bound for $\x_A$ would be left unchanged (red dot). Hence, this ordering would yield a strictly better bound (for the samples $\x_A$ and $\x_B$).  This also implies that the bound with the original ordering is inadmissible.  }
  		\label{fig:XUnionAB}
		 \end{center}
\end{figure}

We now illustrate how the central optimization problem may produce bounds with the same value for two or more consecutive elements of the ordering of a sample space. Unlike in Figure~\ref{fig:incrementalBounds}, where the argmin for the central optimization problem is in the open simplex, we consider the situation illustrated in Figure~\ref{fig:XUnionAB}, where the argmin distribution is on the boundary of the simplex (left panel of Figure~\ref{fig:XUnionAB}). 

In this type of situation, the central optimization problem for two consecutive samples $\x_A$ and $\x_B$ in a total order may produce the same value. This is because the likelihood of a point on the boundary of the simplex (the red point) need not be increased by the addition of a new sample, since the distribution at the red point may assign a likelihood of zero to that new sample. This results in two (or more) consecutive bounds with the same result. 

Notice, however, that if we reverse the order of $\x_A$ and $\x_B$ in the total order to form a new total order in which $\x_B$ precedes $\x_A$, then the sequence of bounds would be governed by the mean of the yellow dot distribution (center of Figure~\ref{fig:XUnionAB}) and then the mean of the red dot, which is lower. That is, by reversing the order of $\x_A$ and $\x_B$ in the total order, we break the tie that was occurring. We now formalize these concepts.

\subsubsection{Breakable and unbreakable ties}
Let $B_T$ be a conditionally optimal bound based on a total order $T$, and suppose that $B_T(\x_{t_k})=B_T(\x_{t_{k+1}})$; that is, the two samples have the same bound value. Let $\Omega_k^T$ and $\Omega_{k+1}^T$ be the upper sets for $\x_{t_k}$ and $\x_{t_{k+1}}$ under the order $T$. Now consider another ordering $U$ (with corresponding bound $B_U$) which is equal to $T$ except that the positions of the elements $t_k$ and $t_{k+1}$ are swapped: $u_k=t_{k+1}$ and $u_{k+1}=t_k$.  And let $\Omega_k^U$ and $\Omega_{k+1}^U$ be the upper sets for $\x_{u_k}$ and $\x_{u_{k+1}}$. 

Note first that $\Omega^T_k=\Omega^U_k$, and so we must have $B_T(\x_{t_k})=B_U(\x_{u_k})$.  We also know that $B_U(\x_{t_k})=B_U(\x_{u_{k+1}})\geq B_U(\x_{u_k})=B_T(\x_{t_k})$. This means that whenever we swap the order of two samples that have tied bounds, under the new bound, one of the bound values will not change, and the other value must improve (get larger) or remain the same.  In other words, for tied bounds, swapping the order of the elements involved in the tie will either lead to an improved bound for one of the samples, or will lead to no change.

Swapping adjacent tied samples, as above, provides one way to break a tie. With three or more tied samples, however, any single adjacent swap in the original order may leave the bound unchanged, while multiple swaps can place a sample in a position whose upper set yields a strictly larger bound. We therefore define breakability using all total orders consistent with the bound. Such orders preserve the relative order of samples assigned unequal bound values but allow tied samples to be reordered freely. They therefore capture all upper sets through which a tie could be broken.

\begin{definition}[Bounds with breakable and unbreakable ties]
\label{def:unbreakable}
Let $B_T$ be a conditionally optimal bound for a total order $T$. We say that $B_T$ has \emph{unbreakable ties} if, for every total order $T'$ with which $B_T$ is consistent, letting $B_{T'}$ denote the associated conditionally optimal bound, we have $B_T(\x) = B_{T'}(\x)$ for each $\x \in \Omega$. Otherwise, we say that $B_T$ has a \emph{breakable tie}; equivalently, there exists such a $T'$ for which $B_{T'} \neq B_T$.
\end{definition}

Thus, an improvement after a single adjacent swap demonstrates a breakable tie. Invariance under one such swap does not, however, establish unbreakability, which requires invariance under every possible reordering of the samples within tied groups.

\subsubsection{Implications of breakable and unbreakable ties}
Intuitively, a breakable tie allows the tied samples to be reordered in a way that strengthens the bound without weakening it elsewhere. The original bound therefore cannot be admissible, as formalized in the following lemma.

\begin{lemma}[Bounds with breakable ties are not admissible]
\label{lem:breakable}
Let $B_T$ be the conditionally optimal bound for a total order $T$. If $B_T$ has a breakable tie, then it is not admissible.
\end{lemma}
\begin{proof}
    By Definition~\ref{def:unbreakable}, there exists a total order $T'$, with which $B_T$ is consistent, such that the associated conditionally optimal bound $B_{T'}$ differs from $B_T$. Since $B_T$ is consistent with $T'$, Theorem~\ref{thm:conditional_optimality} implies that $B_{T'} \geq B_T$. And because $B_{T'} \neq B_T$, there must exist at least one sample $\x \in \Omega$ for which $B_{T'}(\x) > B_T(\x)$. It follows then that $B_{T'}$ dominates $B_T$. Hence $B_T$ is not admissible.
\end{proof}

\begin{lemma}[Bounds with unbreakable ties are admissible]
\label{lem:unbreakable}
    Let $T$ be a total order and $B_T$ the conditionally optimal bound associated with $T$. If $B_T$ has unbreakable ties, then it is admissible with respect to the full family of lower bound functions over $\Omega$.
\end{lemma}

\begin{proof}
    Suppose for the sake of contradiction that some bound $C$, consistent with a total order $U$, dominates $B_T$. Then by Theorem~\ref{thm:conditional_optimality} we must have that $B_U$, the conditionally optimal bound consistent with $U$, also dominates $B_T$. Since $B_T$ has unbreakable ties by assumption, we have that $B_U  = B_T$ whenever $B_T$ is consistent with $U$, which contradicts the dominance of $B_U$. Thus, we proceed under the assumption that $B_T$ is not consistent with $U$.

    We can construct a new total order $T'$ with which $B_T$ is consistent as follows. $T'$ maintains the strict order induced by $B_T$, i.e., $\forall \x, \y \in \Omega$, if $B_T(\x) < B_T(\y)$, then $\x$ precedes $\y$ in $T'$. And for any block of samples that are tied according to $B_T$, $T'$ follows the order $U$. Let $B_{T'}$ be the conditionally optimal bound associated with $T'$.  Because $B_T$ has unbreakable ties by assumption, we have that $B_T = B_{T'}$.
    
    For orders $U$ and $T'$, label the samples of $\Omega$ as follows.
    \begin{eqnarray}
        U = (\y_1, \ldots, \y_N), \qquad T' = (\z_1, \ldots, \z_N).
    \end{eqnarray}
    Beginning with the lowest samples in each order, let $k$ be the first index where $y_k \neq z_k$. Samples that are tied according to $B_T$ follow the order $U$; therefore, a deviation from order $U$ at index $k$ implies that $B_T(\z_k) < B_T(\y_k)$. Of course we also have that $B_{T'}(\z_k) = B_T(\z_k)$ by the unbreakability of ties in $B_T$ and the fact that $B_T$ is consistent with $T'$ by construction. But because $U$ and $T'$ agree on the first $k-1$ samples, their upper sets for the $k$th sample are identical. It follows then that
    \begin{eqnarray}
        B_U(\y_k) = B_{T'}(\z_k) = B_T(\z_k) < B_T(\y_k),
    \end{eqnarray}
    which contradicts the dominance of $B_U$ over $B_T$. Therefore, no valid bound dominates $B_T$, making it admissible.
\end{proof}

\begin{theorem}
\label{thm:admissibility_characterization}
    Let $B^*_T$ be the conditionally optimal bound for total order $T$. Then $B^*_T$ is admissible if and only if for any other total order $U$ with which $B_T^*$ is consistent and whose conditionally optimal bound is $B_U^*$,
    \begin{eqnarray}
        B^*_T = B^*_U.
    \end{eqnarray}
\end{theorem}

\begin{proof}
    Suppose first that $B^*_T$ is admissible but $B^*_T \neq B^*_U$. Since $B^*_T$ is consistent with $U$, and $B^*_U$ is conditionally optimal by assumption, we must have that $B^*_U \geq B^*_T$. But since $B^*_T \neq B^*_U$, it must be the case that $B^*_U$ dominates $B^*_T$, contradicting admissibility. It follows then that $B^*_T = B^*_U$.
    
    Next, suppose that $B^*_T = B^*_U$ for each $U$ with which $B^*_T$ is consistent. Then $B^*_T$ has unbreakable ties according to Definition~\ref{def:unbreakable}. Thus, it must also be admissible by Lemma~\ref{lem:unbreakable}.
\end{proof}
\color{black}

Summarizing the results so far, we can say that all conditionally optimal injective bounds are admissible (Lemma~\ref{lem:injective_optimal_bounds}). Among non-injective conditionally optimal bounds, those with a breakable tie are not admissible (Lemma~\ref{lem:breakable}), whereas those with unbreakable ties are admissible (Lemma~\ref{lem:unbreakable}).

This completes our characterization of admissibility. These results of course lead to new questions such as which total orders lead to admissible bounds and which do not.  One perhaps counterintuitive result that we shall not dive into here is that there is no need for a total order to obey the ``natural partial order'' on samples in order for it to be admissible. By natural partial order, we mean the coordinatewise order defined by $\x\leq\y \iff x_i\leq y_i \; \forall i$. There are many total orderings which do not obey such a natural partial order, and nevertheless lead to admissible bounds. We leave this as a topic for future work.

However, we will discuss below a class of total orderings that {\em never} leads to an admissible bound. 
In the next subsection, we consider the special case of conditionally optimal bounds that depend upon an ordering in which the lowest possible sample is {\em not} the first sample in the total order over the sample space. We refer to these as {\em degenerate} bounds, and they illustrate many of the concepts discussed above.

\subsubsection{The special sample (0,0,...,0)}
For simplicity of exposition, and without loss of generality, let us assume that the least element of the support $S$ is $0$. We define the unique sample $\mathbf{0}=(0,0,...,0)$ of a sample space as the sample all of whose components are 0. As we shall see, bounds based on a total order which does not put $\mathbf{0}$ as the first element have a certain degenerate behavior.

\begin{definition}[Degenerate total order and degenerate bound]
Consider a conditionally optimal lower bound over a sample space $\Omega$ with support set $S$ with least element $0$, that is specified by a total order $T$ on the sample space, with first element $\x_{t_1}$. If $\x_{t_1} \neq \mathbf{0}$, we say that the total order and the resulting bound are ${\bf degenerate}$; if $\x_{t_1}=\mathbf{0}$, we say that the total order and the resulting bound are {\bf non-degenerate}.  
\end{definition}

\begin{definition}[vacuous bound]
    For a given sample $\x$, we say that a lower bound is {\bf vacuous} if its value is equal to the minimum of the support (here, assumed to be $0$). If the minimum of the support is $0$, it is already known that $\mu\geq 0$ before a sample is seen, so a vacuous bound provides no new information.
\end{definition}

\begin{lemma}[Bounds for probabilities of subsets of $\Omega$ containing $\mathbf{0}$] 
\label{lemma:leastSample}
 Consider the unique sample $\mathbf{0}=(0,0,...,0)$ of a sample space $\Omega$. Let $\Omega_{\mathbf{0}}$ be a subset of $\Omega$ that contains $\mathbf{0}$. And let $\mathcal{F}(\Omega_{\mathbf{0}},\alpha)= cl(\{F\in \mathcal{F}: \Pr_F(\Omega_{\mathbf{0}}) > \alpha\})$. For any such $\Omega_{\mathbf{0}}$, we have that
\begin{eqnarray}
\min_{F \in \mathcal{F}(\Omega_{\mathbf{0}},\alpha)} \mu(F)=0.
\end{eqnarray}
\end{lemma}
\begin{proof}
Let $F_0 \in \mathcal{F}$ be the distribution with all of its mass on $0$. Since 
$$
\Pr\nolimits_{F_0}(\mathbf{0}) = 1,
$$
and
$\Omega_{\mathbf{0}}$ contains $\mathbf{0}$, we have
that 
$$
\sum_{\x\in \Omega_{\mathbf{0}}} \Pr\nolimits_{F_0}(\x) =1,
$$
which is greater than or equal to $\alpha$ for any value of $\alpha$. 
We also have that $\mu(F_0) = 0$, which is the lowest possible mean for \emph{any} distribution in the simplex. Therefore, for any subset $\Omega_{\mathbf{0}}$ containing $\mathbf{0}$, the lower bound will be 0.
\end{proof}

\begin{corollary}
\label{cor:vacuous}
    If $\mathbf{0}$ is the $k$th element of a total order $T$, then the degenerate bound based on this total order will have vacuous bounds (of 0) for the first $k$  samples in the sample space.
\end{corollary}
\begin{proof}
    According to Lemma~\ref{lemma:leastSample}, any subset of a sample space that contains the sample $\mathbf{0}$ yields a bound of $0$ (as a result of the central optimization problem). If $\mathbf{0}$ is the $k$th element of a total order, then the first $k$ subsets $\Omega_i$ used in Equation~\ref{eq:CentralOptimizationProblem} will all contain $\mathbf{0}$, leading to bounds of 0 for the first $k$ elements of the total order $T$. 
\end{proof}

\begin{lemma}[non-degenerate, conditionally optimal bounds are vacuous only for the sample $\mathbf{0}$]
\label{lemma:non-degenerate}
A non-degenerate, conditionally optimal bound takes the value 0 at its lowest-ranked sample, as established by Lemma~\ref{lemma:leastSample}. That is, its lowest bound is vacuous. However, none of the other bounds can be vacuous. That is, they are all greater than 0.
\end{lemma}
\begin{proof}
Without loss of generality, we assume that the least element of the support of the family of distributions is 0. To prove this lemma, we consider under what circumstances the central optimization problem produces 0. That is, for what sample spaces $\Omega$, subsets $\Phi$ of $\Omega$, and confidence levels $1-\alpha$  can
\begin{eqnarray}
\label{eq:inf_eq_0}
\underset{cl\{F \in \mathcal{F}: \Pr_F(\Phi) > \alpha\}}{\min} 
\mu(F)=0
\end{eqnarray}
hold true?
Since we are only considering non-degenerate bounds, we assume that the sample $\bf{0}$ is first in the ordering and thus not in the upper set $\Phi$ for any other sample $\x\neq \mathbf{0}$. 

Consider a specific upper set $\Phi$ such that $\mathbf{0}\notin \Phi$. Let
$$\mathcal{F}_\Phi=cl\{F \in \mathcal{F}: \Pr\nolimits_F(\Phi) > \alpha\}.$$
Note that every distribution in $\mathcal{F}_{\Phi}$ assigns probability at least $\alpha$ to $\Phi$. 

We start by observing that $F_0$ is the only distribution in  $\mathcal{F}$ whose mean is $0$.  Thus, Equation~\ref{eq:inf_eq_0} can only hold true if $F_0\in \mathcal{F}_{\Phi}$.

Suppose $F_0\in \mathcal{F}_\Phi.$ For any $\alpha>0$, the probability of at least one sample in $\Phi$ must be non-zero to have $\Pr\nolimits_{F_0}(\Phi)\geq \alpha$. However, we have
$$
\Pr\nolimits_{F_0}(\x) = 0 \iff \x\neq \mathbf{0}.
$$
That is, $F_0$ assigns a zero probability to {\em every sample} except the zero sample itself. And since by assumption $\mathbf{0} \notin \Phi$, $\Pr_{F_0}(\Phi)=0.$
Thus $F_0$ does not meet the criterion to be in $\mathcal{F}_\Phi$, contradicting the assumption.
\end{proof}

\begin{lemma}
\label{lem:degen_are_dom}
Given a degenerate conditionally optimal bound based on a total order $T$, there exists another conditionally optimal bound, based on another order $T'$, that dominates it. Hence, no degenerate bound is admissible. 
\end{lemma}
\begin{proof}
Given a degenerate conditionally optimal bound $A$ based upon a total order $T$, there exists a bound that dominates $A$ based on a modified total order $T'$. Let $k\neq 1$ be the index of the sample $\bf{0}$ in the order $T$. (Recall that the index of $\mathbf{0}$ is $1$ only for non-degenerate bounds.) Corollary~\ref{cor:vacuous} tells us that the first $k$ bound values of $A$ will be $0$. If we create a new total order $T'$ by swapping the positions of $t_1$ and $t_k$, then a bound $B$ conditioned on the new order will be non-degenerate. According to Lemma~\ref{lemma:non-degenerate}, the bound values for samples $2$ through $k$ of $B$ will be non-zero, and hence are stronger than the bounds for $A$. In addition, since the upper sets of $T$ and $T'$ are equivalent for positions $k+1$ and later in these orderings, their bounds will be equivalent. Hence, $B$ dominates $A$ for samples $2$ through $k$ and is equivalent for the rest, meaning that the bound $B$ dominates the bound $A$.  
\end{proof}

In the language of breakable and unbreakable ties, any degenerate bound has a breakable tie and hence cannot be admissible.

\begin{corollary}(upper bound on number of admissible bounds)
\label{cor:max_admis}
If $N$ is the size of a sample space $\Omega$, there are no more than $(N-1)!$ admissible bounds over $\Omega$.
\end{corollary}
\begin{proof}
We can tighten the previous upper bound on the number of admissible bounds. Since total orders must begin with the $\mathbf{0}$ sample in order to be admissible, this leaves us with only $(N-1)!$ possible orderings that might satisfy the admissibility requirement. 
\end{proof}

\section{On the non-existence of optimal bounds}

In this section, we address the question of whether there exists an optimal bound, i.e., a bound that dominates all other valid bounds, for a given sample space and confidence level. We show that, except in trivial cases, there is no valid lower bound that dominates every other valid lower bound. 
To determine whether an optimal bound exists, it suffices, by Lemma~\ref{lem:admissible_bounds} and Corollary~\ref{cor:nfactorial}, to consider the finite set of conditionally optimal bounds. The key observation is that every such bound lies on a finite dominance chain ending at an admissible bound. In particular, if the current bound is not itself admissible, we may replace it with another conditionally optimal bound that dominates it. Because there are finitely many such bounds and each replacement is a strict improvement, this process must eventually reach an admissible bound. This observation allows us to rule out the existence of an optimal bound by starting from two conveniently chosen conditionally optimal bounds and showing that their dominance chains terminate at distinct admissible bounds.

\subsection{Conditions for the non-existence of optimal bounds}

By Definition~\ref{def:opt}, if there are two or more admissible bounds over a sample space $\Omega$, then there is no optimal bound. Below we show that this is the case for all but trivial finite sample spaces and choices of $\alpha$. 
The same conclusion can be drawn from \cite{KabailaLloyd2000} (Example 3) for $\alpha \in (0, 1/2)$. Here we establish it for every $\alpha \in (0, 1)$ by explicitly constructing two distinct admissible bounds for each nontrivial finite sample space.

In this section we make use of two canonical samples $\x_A \coloneqq (1, \ldots, 1)$ and $\x_B \coloneqq (0, 1, \ldots, 1)$. We begin with a lemma establishing that $\x_A$ and $\x_B$ create strictly increasing bounds when they are placed consecutively at the end of an order.

\color{black}
\begin{lemma}
\label{lem:consec_samples}
    If $A$ and $B$ are two conditionally optimal bounds such that $A$ places $\x_A$ last and $\x_B$ second to last, and $B$ places $\x_B$ last and $\x_A$ second to last, then $B(\x_A) < A(\x_A)$ and $A(\x_B) < B(\x_B)$.
\end{lemma}

\begin{proof}
    We begin by observing that, for any $F \in \mathcal{F}$, moving all mass that $F$ assigns to $S \setminus \{0,1\}$ to 0 cannot increase its mean, leaves $\Pr_F(\x_A)$ unchanged, and cannot decrease $\Pr_F(\x_B)$. Consequently, the infima defining the Buehler bounds for the three relevant upper sets $\{\x_A\}$, $\{\x_B\}$, and $\{\x_A,\x_B\}$ are unchanged when the optimization is restricted to distributions supported on $\{0,1\}$. Thus, to show that one bound is strictly greater than another for any given sample, it will suffice by Lemma~\ref{lem:suff_cond_no_tie} to show that the distribution achieving the bound's value will be in the one-dimensional open simplex. 
    
    To that end, let $p = \Pr_F(1) = \mu(F)$. We have that 
    \begin{eqnarray}
    A(x_A)&=&\min_{cl\{F: \Pr\nolimits_F(\x_A)> \alpha\}} \;p\\
    &=& \min_{cl\{F: p^n> \alpha\}} \;p \\
    &=& \min_{cl\{F: p> \alpha^{1/n}\}}\; p \\
    &=& \alpha^{1/n}.
    \end{eqnarray}
    Note that for all values of $\alpha$ between $0$ and $1$, we have that $\alpha^{1/n}$ is also between $0$ and $1$, so the distribution achieving this minimum will be in the one-dimensional open simplex. It must then be the case that $A(\x_A)>B(\x_A)$.

    We argue similarly for the bound $B$. Again let $p = \Pr_F(1)$. Suppose first that there exists no $p \in [0,1]$ such that $\Pr_F(\x_B)> \alpha$. Then $B(\x_B) = \infty$ so that clearly $A(\x_B) < B(\x_B)$. If such a $p$ does exist, then we have
    \begin{eqnarray}
    B(\x_B)&=&\min_{cl\{F: \Pr\nolimits_F(\x_B)> \alpha\}} \; p\\
    &=& \min_{cl\{F: n(1-p)p^{n-1} > \alpha\}} \; p
    \end{eqnarray}
    $B(\x_B)$ is the infimum of values $p$ such that $n p^{n-1} (1-p) > \alpha$. And since $\alpha \in (0,1)$, that infimum will be achieved neither at $p=0$ nor $p=1$. Thus, the resulting distributions of this optimization, when they exist, must also be in the open simplex. Therefore, $A(\x_B) < B(\x_B)$.
\end{proof}

Having established the properties of bounds placing $\x_A$ and $\x_B$ at the end of their orders, we now use them to construct two distinct admissible bounds on the same sample space. Starting from each bound in turn, we repeatedly replace the current bound by a conditionally optimal one that dominates it until an admissible bound is reached.

\color{black}
\begin{theorem}[Non-existence of optimal bounds]
\label{thm:no_opt}
Let $\Omega(S,n)$ be a sample space with a support set $S$ of at least 2 elements, and let $n$ be a sample size of at least 2. Then there exists no optimal lower bound for $\Omega(S,n)$ for any confidence $0<1-\alpha<1$.
\end{theorem}

\begin{proof}
    Any sample space is simply an affine transformation of a sample space whose extremal elements are 0 and 1. Therefore, it suffices to prove the conclusion subject to this restriction. To that end, we restrict our attention to $S$ where $\min S = 0$ and $\max S = 1$. Let $A$ and $B$ be bounds as defined in Lemma~\ref{lem:consec_samples}. Corresponding to $A$ is a sequence of conditionally optimal bounds such that 
    $$
    A \leq \ldots \leq A_{\max};
    $$
    a similar sequence exists for $B$. By construction, neither $A_{\max}$ nor $B_{\max}$ is dominated by any other bound. Thus, to prove that there exists no optimal bound, it will suffice by Definition~\ref{def:opt} to show that $A_{\max} \neq B_{\max}$.

    Since $A_{\max} \geq A$, we have that $A_{\max}(\x_A) \geq A(\x_A)$ (due to domination), but also $A_{\max}(\x_A) \leq A(\x_A)$, since $A(\x_A)$ is the highest possible bound for $\x_A$, so $A_{\max}(\x_A)=A(\x_A)$. Now suppose that $A_{\max}(\x_B) > A(\x_B)$. Since we have already argued that $A_{\max}(\x_A)=A(\x_A)$, Lemma~\ref{lem:consec_samples} establishes that $A_{\max}(\x_A) > A(\x_B)$ as well. Let $\x'$ be whichever sample among $\x_A$ and $\x_B$ appears lower in the order induced by $A_{\max}$, and take $\Omega'$ to be its upper set according to the same order. We have already established that $A_{\max}(\x') > A(\x_B)$. Since both $A$ and $A_{\max}$ are valid, Theorem~\ref{thm:conditional_optimality} implies that $A(\x_B)$ is the Buehler bound for $\{\x_A, \x_B\}$ and $A_{\max}(\x')$ is the Buehler bound for $\Omega'$.
    Although these upper sets arise from different total orders, the value of their respective Buehler bounds depends only on the upper sets themselves, so the argument of Lemma~\ref{lemma:incremental_sums} applies: since $\{\x_A, \x_B\} \subseteq \Omega'$ we must have $A_{\max}(\x') \leq A(\x_B)$, a contradiction. Thus we must have $A_{\max}(\x_B) \leq A(\x_B)$, and $A_{\max}(\x_B) = A(\x_B)$ follows by the fact that $A_{\max}$ dominates $A$.
    
    A similar argument shows that $B_{\max}(\x_A) = B(\x_A)$ and $B_{\max}(\x_B) = B(\x_B)$. Because $A(\x_A) > B(\x_A)$, again by Lemma~\ref{lem:consec_samples}, it immediately follows that $A_{\max}(\x_A) > B_{\max}(\x_A)$, which proves that $A_{\max} \neq B_{\max}$. Therefore, there exists no optimal bound.
\end{proof}
\color{black}

\section{Summary of results}

We have studied valid lower confidence bounds for the mean that are based on $n$ i.i.d. observations from an unknown distribution  supported on a known finite set $S \subset \mathbb{R}$. Our results build fundamentally on the samplewise dominance relationship between bounds: one bound dominates another if it is at least as large for every possible sample and strictly larger for at least one. Conditioning on a total order of the $n$-trial multinomial sample space over categories $S$ reduces the construction of a bound to a sequence of nested upper-set optimization problems. For each order, the classical Buehler construction takes the infimum mean over distributions assigning probability exceeding $\alpha$ to the corresponding upper set and yields the greatest valid bound consistent with that order. Our main contribution is to determine which of these fixed-order optima remain admissible when compared with all valid bounds.

We showed that every admissible bound must be conditionally optimal for at least one total order. For such a bound $B_T^*$, global admissibility holds if and only if every total order $U$, disagreeing with order $T$ only among samples tied in $B_T^*$, produces the same conditionally optimal bound, i.e., $B_U^* = B_T^*$. Consequently, every injective, conditionally optimal bound is admissible, whereas a breakable tie guarantees inadmissibility. Another interpretation is that a non-injective, conditionally optimal bound is admissible precisely when it is invariant under every compatible reordering of its tied groups. This provides a complete structural characterization of the admissible frontier under samplewise dominance.

Our work also constrains the number of admissible bounds over a given sample space.  Any order producing an admissible bound must begin with the sample consisting entirely of the minimum support value. Thus, if $N$ is the size of the sample space, then there are at most $(N-1)!$ distinct admissible bounds for each fixed $\alpha \in (0,1)$. Nevertheless, the admissible frontier is typically not a singleton. Whenever $|S| \ge 2$, $n \ge 2$, and $\alpha \in (0,1)$, at least two distinct admissible bounds exist. Hence no valid lower confidence bound samplewise dominates every other valid lower confidence bound in this nontrivial regime. 

Several questions remain. Although the characterization is finite, the number of candidate orders grows factorially, and evaluating an individual candidate may require a difficult optimization over the probability simplex. Important directions include identifying orders whose conditionally optimal bounds are admissible and tractable to compute, developing approximations that retain bound validity, and extending the present theory to infinite or uncountable support.
\color{black}

\section{Disclosure}

AI tools were used in preparing this document to support literature searches, clarify the relationship between prior work and the paper’s contributions, refine prose, and assist with corrections to existing proofs. The authors take full responsibility for the content, including any errors or omissions.

\bibliographystyle{plainnat}  
\bibliography{mean_interval}
\end{document}